\documentclass[10pt]{article}
\usepackage{amsfonts,amsmath,amssymb,color}
\usepackage{ulem}
\usepackage[top=.95 in, bottom = 1 in, left=1in, right = .6in]{geometry}
\usepackage{anysize}
\usepackage[title]{appendix}

\renewcommand{\em}{\it}
\newcommand{\R}{{\mathbb R}}
\newcommand{\E}{{\mathbb E}}
\newcommand{\B}{{\mathbb B}}
\renewcommand{\P}{{\mathbb P }}

\newcommand{\N}{{\mathbb N }}
\newcommand{\PA}{{\mathcal P}}
\newcommand{\G}{{\mathcal G}}
\newcommand{\LA}{{\mathcal L}}
\newcommand{\FA}{{\mathcal F}}
\newcommand{\1}{\textbf{1}}
\newcommand {\ben}{\begin{equation}}
\newcommand {\een}{\end{equation}}
\newcommand {\bena}{\begin{eqnarray}}
\newcommand {\eena}{\end{eqnarray}}
\newcommand {\benas}{\begin{eqnarray*}}
\newcommand {\eenas}{\end{eqnarray*}}

\newcommand{\norm}[1]{\parallel#1\parallel}

\newtheorem{theorem}{Theorem}[section]
\newtheorem{lemma}[theorem]{Lemma}

\newtheorem{proposition}[theorem]{Proposition}
\newtheorem{example}{Example}
\newtheorem{remark}[theorem]{Remark}

\newtheorem{definition}[theorem]{Definition}

\def\qed{\hfill $\Box$}

\definecolor{nb}{rgb}{.6,.176,1}
\definecolor{red}{rgb}{.8,0,0}
\definecolor{darkgreen}{rgb}{0,.5,0}
\definecolor{Red}{rgb}{1,0,0}
\definecolor{Blue}{rgb}{0,0,1}
\definecolor{Olive}{rgb}{0.41,0.55,0.13}
\definecolor{Yarok}{rgb}{0,0.5,0}
\definecolor{Green}{rgb}{0,1,0}
\definecolor{MGreen}{rgb}{0,0.8,0}
\definecolor{DGreen}{rgb}{0,0.55,0}
\definecolor{Yellow}{rgb}{1,1,0}
\definecolor{Cyan}{rgb}{0,1,1}
\definecolor{Magenta}{rgb}{1,0,1}
\definecolor{Orange}{rgb}{1,.5,0}
\definecolor{Violet}{rgb}{.5,0,.5}
\definecolor{Purple}{rgb}{.75,0,.25}
\definecolor{Brown}{rgb}{.75,.5,.25}
\definecolor{Grey}{rgb}{.7,.7,.7}
\definecolor{Black}{rgb}{0,0,0}

\definecolor{violet}{rgb}{.3,.0,.9}

\begin{document}

\title{Simultaneous Small Noise Limit for Singularly Perturbed Slow-Fast Coupled Diffusions.}
\author{Siva R. Athreya, Vivek S. Borkar, K. Suresh Kumar, and Rajesh Sundaresan
}
\maketitle

\begin{abstract}
We consider a simultaneous small noise limit for a
singularly perturbed coupled diffusion described by
\begin{eqnarray*}
dX^{\varepsilon}_t &=& b(X^{\varepsilon}_t, Y^{\varepsilon}_t)dt + \varepsilon^{\alpha}dB_t,
\\
 dY^{\varepsilon}_t &=&  -  \frac{1}{\varepsilon}  \nabla_yU(X^{\varepsilon}_t, Y^{\varepsilon}_t)dt +
\frac{s(\varepsilon)}{\sqrt{\varepsilon}} dW_t,
\end{eqnarray*}
where $B_t, W_t$ are independent Brownian motions  on $\R^d$ and $\R^m$ respectively, $b : \mathbb{R}^d
\times \mathbb{R}^m \rightarrow \mathbb{R}^d$, $U : \mathbb{R}^d
\times \mathbb{R}^m \rightarrow \mathbb{R}$ and $s :(0,\infty)
\rightarrow (0,\infty)$.  We impose regularity assumptions on $b$, $U$ and let  $0 < \alpha < 1.$ When $s(\varepsilon)$ goes to zero slower than a prescribed rate as  $\varepsilon \rightarrow 0$, we characterize all weak limit points of $X^{\varepsilon}$, as $\varepsilon \rightarrow 0$, as
solutions  to a differential equation driven by a measurable vector field. Under an additional   assumption on the behaviour of $U(x, \cdot)$ at its global minima we   characterize all limit points as Filippov solutions to the  differential equation.
\end{abstract}


\section{Introduction}
In this article we consider the simultaneous small noise limit for a
singularly perturbed coupled slow-fast diffusion  given by
\begin{eqnarray}
dX^{\varepsilon}_t &=& b(X^{\varepsilon}_t, Y^{\varepsilon}_t)dt + \varepsilon^{\alpha}dB_t,
\label{ex0}\\
 dY^{\varepsilon}_t &=&  -  \frac{1}{\varepsilon}  \nabla_yU(X^{\varepsilon}_t, Y^{\varepsilon}_t)dt +
\frac{s(\varepsilon)}{\sqrt{\varepsilon}} dW_t, \label{wye0}
\end{eqnarray}
where $B_t, W_t$ are independent Brownian motions on $\R^d$ and $\R^m$ respectively,  $b : \mathbb{R}^d
\times \mathbb{R}^m \rightarrow \mathbb{R}^d$, $U : \mathbb{R}^d
\times \mathbb{R}^m \rightarrow \mathbb{R}$ and $s :(0,\infty)
\rightarrow (0,\infty)$. We impose regularity assumptions on $b$,$U$ and
let $0 < \alpha < 1$. When $s(\varepsilon)$ goes to zero
slower than a prescribed rate as $\varepsilon \rightarrow 0$, we show that in the simultaneous small-noise limit all weak limit points $X$ satisfy
\begin{equation} \label{odemain0}
\frac{d}{dt}X_t = \int_{\R^d} b(X_t, y)\nu^{0,X_t}_t(dy)
\end{equation}
where $\nu^{0,X_t}_t(dy)$ is a probability measure supported on
finitely many global minima of $U(X_t, \cdot)$ (see Theorem \ref{maintheorem1}). If an additional assumption on the behaviour
  of $U(x,\cdot)$ at its global minima is made then we show that
  $\nu^{0,X_t}_t(dy)$ is time independent and given by a determinantal
  formula arising from Laplace's principle. Consequently, for this
  class of $U$, we show that every limit point is a generalized
  Filippov solution to (\ref{odemain0}) driven by a vector field (see
  Theorem \ref{maintheorem2}).  In Section \ref{mr} we state the model,
  assumptions made, and the two
  results  precisely and in Section \ref{examples-refinements} we discuss
  examples of $U$ that satisfy the required assumptions.

The factor $\frac{1}{\varepsilon}$ in the drift term in
\eqref{wye0}, intuitively suggests that the $Y^\varepsilon$ process is
the ``fast moving'' process as $\varepsilon \rightarrow 0$ and that
the ``slow moving'' process $X^\varepsilon$ will see an averaging of
$Y$ in this limit. The study of averaging principle in various
dynamical systems dates back to the work of Khasminskii and others, summarized in, e.g.,  Freidlin and Wentzell
\cite{FW12}, Kabanov and Pergamenshchikov \cite{KP03}. The dynamical systems considered there involve a ``slow
process'' $X^{\varepsilon}$ as a solution to an ordinary differential
equation (i.e.~\eqref{ex0} with no $B_t$ term) coupled with the fast
process $Y^{\varepsilon}$ given by a stochastic differential equation
with absence of small noise (i.e.~\eqref{wye0} with $s(\varepsilon) =1
$).  In this setting, under further assumptions on $b,U$, the
averaging principle leading to characterization of limit points,
normal deviations, and large deviations from the averaging principle
are detailed in \cite[Chapter~7]{FW12}. The ground work for this lies
in understanding the long-term behavior of solutions to \eqref{wye0}
(for fixed $\varepsilon >0$), and is laid out in
\cite[Chapters~4-6]{FW12}. We shall rely on this foundation in
prescribing assumptions for $U$ in our main results.

Large deviations and generalizations to ``full dependence'' systems
were considered in the works of Veretennikov in \cite{V13, V94,
  V99}. Motivated by questions from homogenization, \cite{V00}
considered the fast process \eqref{wye0} with $s(\varepsilon) = 1$ but with presence of
small noise for the slow process (i.e.\  \eqref{ex0} with $\alpha=\frac{1}{2}$) and established a large deviation principle (LDP) for
$X^{\varepsilon}$ as $\varepsilon \rightarrow 0$. One can characterize
the limit points of $X^{\varepsilon}$ as $\varepsilon \rightarrow 0$
as the set where the rate function is equal to zero (see \cite[Remark
  3]{V00}). In \cite{L96}, Liptser considered the joint distribution of
the slow process and of the empirical process associated with the fast
variable in the one-dimensional setting and derived an LDP. This was
recently generalized to multidimensional and full dependence systems by Puhalskii in \cite{P14}.
The diffusions driving the slow and the fast processes in \cite{P14} do not have to be uncorrelated.

In related works, Spiliopoulos in \cite{S13,S14}, Morse and Spiliopolous in \cite{MS17}, and Gailus and Spiliopoulous in \cite{GS17} considered a class of coupled diffusions with multiple time scales in the full dependence setting. Contained therein, after suitable relabelling of the parameters and appropriate choice of coefficients, are results that will apply to (\ref{ex0})-(\ref{wye0}) for specific $b, \nabla_y U$ and with $s(\varepsilon) = \varepsilon^{\alpha -\frac{1}{2}}$.  Thus, when $\alpha < \frac{1}{2}$ the fast process then undergoes stochastic homogenization (i.e. \eqref{wye0} with $s(\varepsilon) \rightarrow \infty$ as $\varepsilon \rightarrow 0$), when $\alpha > \frac{1}{2}$ the fast process has a small noise limit (i.e. \eqref{wye0} with $s(\varepsilon) \rightarrow 0$ as $\varepsilon \rightarrow 0$) and when $\alpha = \frac{1}{2}$ this corresponds to  $s(\varepsilon) = 1$ in \eqref{wye0}. In \cite{S13}, an LDP is shown for the slow process under periodicity assumptions for all the three regimes. Without the periodicity assumption on the coefficients, in \cite{S14} fluctuation results for the slow process are shown in the homogenization and $s(\varepsilon) =1$ regimes, while in \cite{MS17} moderate deviations for the slow process are shown for these two regimes. In \cite{GS17}  parameter estimation results are obtained when $s(\varepsilon) =1$.


Our model falls in the complement of the above.  To the best of our
knowledge the case where no periodicity assumptions are made and when
both slow and fast motions are subjected to small noise limits (i.e.
$\alpha >0$ and $s(\varepsilon) \rightarrow 0$) has not been studied
in the literature. In this paper, we provide a first step towards
understanding this regime. Since we do not impose any periodicity
assumptions on the coefficients this does not allow us to restrict
dynamics on a torus. Thus we have to handle the nontrivial
technicalities that come with a noncompact state space which requires
a new approach.

Our motivation to study this problem comes from a general philosophy
of a selection principle for ill-posed dynamics, attributed to
Kolmogorov in \cite{ER85a}, that adds noise to the dynamics and looks
at the small noise limit for candidate `physical' solution(s). This
philosophy has been variously used in nonlinear circuits \cite{S83},
evolutionary games \cite{FY90}, and underlies the notion of `viscosity
solutions' \cite{FS06}. The problems of `averaging' two time scale
diffusions in the limit of infinite time scale separation on the one
hand \cite{KP03} and of small noise asymptotics for diffusions in the
vanishing noise limit on the other hand \cite{FW12} have been
extensively studied. Our aim here is to analyze the co-occurrence of
the two when the time scale separation and the small noise variance
are controlled by the same parameter $\varepsilon > 0$.

Our first result characterizes any limit point $X$ as a solution to a
differential equation given by (\ref{odemain0}). Inside this result is
contained the interesting observation that the small noise limit in
the faster time scale requires the noise variance to scale in an
inverse logarithmic fashion, or slower (see Remark~\ref{rem:mixing}). In hindsight, this is similar
to the phenomenon observed in optimization algorithms that track the
stationary distribution \cite{CH87}, \cite{GS91}, \cite{HS90} where
the spectral gap determines the convergence rate. So intuitively
speaking not only does the invariant distribution concentrate as the
noise decreases, but also the approach to it slows down because of the
scaling of the second eigenvalue of the infinitesimal generator with
the noise variance. This observation appears to be new under the
additional phenomenon of averaging due to multiple time scales present
in the dynamic itself.

For characterizing the limiting measure in (\ref{odemain0})
   we impose restrictions on the behaviour of $U$ at its global
   minima.  We are then able to identify any limit point as a Filippov
   solution to a differential equation. In particular we are able to
   establish an interesting connection between small noise limits with
   two time scales and the theory of differential equations driven by
   discontinuous vector fields, as in the spirit of single time scale
   case in \cite{BOQ09}. In the single time scale case, there is
   already a considerable body of interesting results, see \cite{BB81, FF14, CH83, BK10}, though a conclusive theory is still wanting.

   We also make an unconventional use of nonlinear filtering theory in
   proving our main result. Nonlinear filtering comes naturally into
   play once we replace the drift of the slow diffusion by its
   conditional expectation given the history of the fast process. It
   is then viewed as the `observation process' in nonlinear filtering
   parlance. We extend the available well-posedness results for
   nonlinear filters to the case when the drift of the `observation'
   process also depends on itself in addition to the `signal'
   process. We prove this in the appendix of this article in Proposition
   \ref{fkksde} and this result is of independent interest (see Remark \ref{rem:filtering}).

We are now ready to state our assumptions and main results in the next subsection.

\subsection{Main Result}
\label{mr}
We use the following notation throughout. For $n \geq 1$, $C_b(\R^n)$ is the space of real valued bounded continuous functions on $\R^n$, $C^2(\R^n)$ is the space of real valued functions with continuous partial derivatives up to
second order, $C^2_b(\R^n) \subset C^2(\R^n)$ are functions in
$C^2(\R^n)$ that in addition are bounded along with their first and second order partial derivatives, and $C^2_0(\R^n) \subset C^2_b(\R^n)$ are functions in
$C^2_b(\R^n)$ that in addition vanish at infinity along with their
first and second order partial derivatives.  We use $\| \cdot \|_2$ for
the $L_2$ norm and  $\| \cdot \|_\infty$ for the sup norm. For a Polish space $S$,
$\PA(S)$ is the Polish space of probability measures on $S$
with the Prohorov topology. For $n \geq 1, x \in \R^n $,  $\norm{x}$ is the usual Euclidean norm, $\langle \cdot , \cdot \rangle$ is the usual inner product, and $\B_1$ is the closed ball of unit radius centered at the origin in that Euclidean space. We use $\nabla_z, D^2_z$ to denote respectively the gradient and the Hessian in variable $z$.

We shall now define the model precisely. Let $0 < \alpha < 1$, $T>0$, $d \geq 1$, $m \geq 1$, $x_0 \in \R^d, y_0 \in \R^m$ be fixed. Let  $(\Omega, \FA, \P)$ be a filtered probability
space  on which $\{B_t\}_{t \geq 0}$ and $\{W_t\}_{t \geq 0}$ are independent standard Brownian motions on $\R^d$ and $\R^m$ respectively. For $0 \leq t \leq T$ and $\varepsilon >0,$ consider the coupled system of stochastic differential equations given by
\begin{eqnarray}
X^{\varepsilon}_t &=& x_0 + \int_0^t b(X^{\varepsilon}_s, Y^{\varepsilon}_s)ds + \varepsilon^{\alpha}B_t,
 \label{ex} \\
 Y^{\varepsilon}_t &=& y_0 -  \frac{1}{\varepsilon} \int_0^t \nabla_yU(X^{\varepsilon}_s, Y^{\varepsilon}_s)ds+
\frac{s(\varepsilon)}{\sqrt{\varepsilon}} W_t,\label{wye}
\end{eqnarray}
where $b : \mathbb{R}^d \times \mathbb{R}^m \rightarrow \mathbb{R}^d$, $U : \mathbb{R}^d \times \mathbb{R}^m \rightarrow \mathbb{R}$, $s : (0, \infty) \rightarrow
(0, \infty)$.

We shall make the following assumptions.
\begin{enumerate}
\item[(B1)] $b \in C_b(\R^d \times \R^m)$  is locally Lipschitz continuous in $y$-variable and is uniformly (w.r.t.\ $y$) Lipschitz continuous in $x$-variable, i.e.   $\exists K_1 > 0  $ such that $\forall \ x, x' \in \R^d, y \in \R^m$
 \begin{equation}\label{ullb}
    \norm{b(x, y) - b(x', y)} \leq K_1 \norm{x-x'}.
  \end{equation}

\item[(U1)]  $U \in C^2(\mathbb{R}^d \times \mathbb{R}^m).$ Further, $\nabla_y U(x, y)$ is  uniformly (w.r.t.\ $y$) Lipschitz continuous in $x$-variable, i.e.  $\exists K_2 > 0 $ such that $\forall \ x, x' \in \R^d$,  $ y \in \R^m$,
  \begin{equation}
    \norm{\nabla_y U(x, y) - \nabla_y U(x', y)}  \leq K_2 \norm{x-x'}. \label{ullu}
  \end{equation}

\item[(U2)]
  There exist $ R >0, M>0, K_3>0 $ such that, for all $x \in \R^d$,
    \bena
    && K_3 \|\xi \|^2 \leq \langle \xi, D_y^2 U(x,y)\xi \rangle  \mbox{ for } \xi \in \R^m \mbox{ and } y \in \R^m, \|y \| > R, \label{U21}\\
    &&\nonumber\\
    && \sup_{\|y \| \leq R} \max\{\mid U(x,y)\mid, \|\nabla_y U(x,y)\|, \|D^2_y U(x,y)\|   \} \leq M, \mbox{ and }\label{U22}\\
        &&\nonumber\\
    && {\ \sup_{y \in \mathbb{R}^m} \Big[ \frac{1}{4\pi e s} \left ( 4 \triangle_y U(x,y)  - \frac{4 \|\nabla_y U (x,y)\|^2}{a} \right ) + 2 U(x,y) \Big]} \leq M s^{\frac{\eta}{\eta-1}}, \nonumber\\
    && \hspace*{2cm} \mbox{ for } a \leq 1, s \geq 1, \mbox{ for some } \eta >1. \label{U23}
    \eena
\end{enumerate}

\begin{remark}  The assumptions (B1) and (U1) immediately imply the local existence and uniqueness of a strong solution for the coupled slow-fast small
diffusions (\ref{ex}) and (\ref{wye}). These along with (\ref{U21}) and
(\ref{U22}) in assumption (U2) imply nonexplosiveness of the system, and thus global existence and uniqueness. Assumption (\ref{U23}) is needed to ensure ultracontractivity (see \cite[Page 363]{BGL14}).

Further, using just (\ref{U21}) and (\ref{U22}) in assumption (U2) we can establish  that  there exists a nonnegative continuous function $g: (0,\infty ) \rightarrow (0,\infty)$ such that
\[
\sup_{z,y \in \R^m: \| z-y \| = r} - \frac{1}{r} \langle \nabla_yU(x,z) - \nabla_yU(x,y), z- y \rangle
  \leq  g(r), \mbox{ for all } r > 0,
\]
with
\begin{equation}
\Gamma:= \int_0^\infty g(s) ds <  \infty.   \label{eqn:PW}
\end{equation}
For completeness, we  provide a proof of  (\ref{eqn:PW}) in  Lemma \ref{l:U24}, Appendix \ref{existenceanduniqueness}. {Along with (\ref{U23}), this is used to obtain a gradient estimate for the fast process.}
\end{remark}

\begin{enumerate}
\item[(U3)]   We assume that $U(x, \cdot)$ has finitely many  critical points for each $x$. For later use, we introduce the following notation for global minima for each $x$ : with $L(x)$ denoting the number of global minima of $U(x, \cdot)$, write
  \bena\label{argminU:x}
  &&  \arg\min U(x, \cdot) :=  \{y_1(x), \cdots , y_{L(x)}(x)\}.
\eena
 Fix $x \in \R^d$. Consider the action functional associated with the ordinary
 differential equation,
\begin{equation*}\label{basicode}
y(t) =  y_0 - \int_0^t  \nabla_y U(x, y(s)) ds, \,
\end{equation*}
defined as follows. For $\varphi \in  C([0, T] ; \mathbb{R}^m)$, write
\begin{equation*}\label{actionfunctional}
  S_T(\varphi) = \begin{cases}  \int^T_0 \norm{\dot{\varphi}(s) + \nabla_y U(x, \varphi(s))}^2 ds &   \mbox{ where } \varphi \mbox{  is absolutely  continuous }\\  & \mbox{ with }  \int^T_0 \norm{\dot{\varphi}(s)}^2 ds < \infty,\\  &\\ \infty & \mbox{ otherwise.}
  \end{cases}
\end{equation*}
Here the dependence of $S_T(\varphi)$ on $x$ is suppressed.
  Define
\begin{equation*}\label{Vfunction}
\tilde{V}(y_i(x), y_j(x)) = \inf \left \{ S_T(\varphi) \left | \begin{array}{l} T>0, \varphi(0) = y_i(x) , \varphi(T) = y_j(x),\\ \varphi(s) \in \mathbb{R}^m \setminus \cup_{k \neq i, j} \{y_k(x)\}, \, 0 \leq s \leq T  \end{array} \right. \right \}.
\end{equation*}

Write $L$ for $L(x)$ and define $[L] := \{1, 2, \ldots , L\}$. For $W \subset [L]$, a graph with node set $[L]$ and directed edges $m \to n$ with
$m \in [L] \setminus W, n \in [L], n \neq m$, is said to be a $W$-graph if
\begin{itemize}
\item each $m \in [L] \setminus W$ is the initial point of exactly one arrow, and
\item there are no cycles in the graph.
\end{itemize}

Let $\G(l), l =1, 2, \ldots, L, $ denote the set of all $W$-graphs with $W$ containing $l$ elements. Set
\begin{equation*}\label{VforMC}
V^{l}(x) \ =  \ \min_{ \chi \in \G(l)} \sum_{(m \to n) \in \chi} \tilde{V}(y_m(x), y_n(x))
\end{equation*}
The additional assumption we require is the following:

\begin{equation} \label{u3}  0 \leq \Lambda:= \sup_{x \in \R^d} \left[ V^1(x) - V^2(x) \right] < \infty.
\end{equation}
\end{enumerate}
  {\begin{remark} Assumption (U3) has two purposes. The first purpose is  as
    in \cite{W72a} and \cite{FW12} to enable the averaging principle
    for the fast process. The second purpose is as in \cite{HS90} to obtain
    spectral gap estimates for speed of convergence of the fast
    process to its invariant measure and to control its rate of equilibration    in the small noise limit via the decay of  $s(\varepsilon)$, see \eqref{spinf}. This brings us to our next assumption.
       \end{remark}}
\begin{enumerate}
\item[(S1)] Our next assumption is on the decay rate of $s(\cdot)$ at $0$. We assume that
 \begin{equation}\label{svarepsilon}
 s(\varepsilon)  \geq  {\sqrt{\frac{C}{\ln(1+ \frac{1}{\varepsilon})}}} \,\, \mbox{ with } \,\,C > \frac{2(\Lambda + 2\Gamma)}{1-\alpha} \,\,  \mbox{ and } \,\,\lim_{\varepsilon \rightarrow 0}  s(\varepsilon) = 0.
  \end{equation}

\end{enumerate}

We are now ready to state the first of the two main results. Recall $C$ from
(\ref{svarepsilon}), $\Gamma$ from (\ref{eqn:PW}) and $\Lambda$ from
(\ref{u3}).

\begin{theorem}\label{maintheorem1}
  Assume {\rm (B1), (U1), (U2), (U3)} and {\rm (S1)}.Then for any
  sequence $\varepsilon_n \downarrow 0$ there is a further
  subsequence, $\varepsilon_{n_k} \downarrow 0$, along which
  $\{X^{\varepsilon_{n_k}}_t, {0 \leq t \leq T}\}$ converges in law on $C([0,T];\R^d)$ to
   $\{X_t, {0 \leq t \leq T}\}$ which is almost surely a
   solution to
  \begin{equation}
X_t = x_0 + \int_0^t \int b(X_s,y) \nu^{0,X_s}_s(dy)ds  \label{odemain1}
\end{equation}
where $\nu^{0,X_s}_s(dy)$ is a probability measure supported on $\arg\min U(X_s, \cdot).$
\end{theorem}

\begin{remark} \label{caratheodory}
  In the proof of the result we shall show that the mapping $s \rightarrow \int
  b(X_s,y)\nu_s^{0,X_s}(dy)$ is almost surely uniformly bounded and  integrable. Thus $X$ is also a Carath\'{e}odory solution\footnote{Carath\'{e}odory solutions
relax the classical requirement that the solution must follow the direction of the vector field at all
times: the differential equation need not be satisfied on a set of measure zero on $[0,T]$. See \cite{BS96} for a precise definition.} to
\begin{align*}
\frac{d}{dt}X_t = \int  b(X_t,y)\nu_t^{0,X_t}(dy),
\end{align*}
with $ X_0 = x_0.$
\end{remark}

We note that the measure $\nu^{0,X_t}_t(dy)$ in the above result may in general
depend on the subsequential limit that is taken. A complete characterization of $\nu_t^{0,X_t}$ is possible in some
 special cases using Laplace's method. For this we impose the following
 additional assumption on the behaviour of $U$ at its global minima.
\begin{enumerate}
\item[(U4)] {For $i \geq 1,$ let
$$D_i = \{ x \in \R^d : L(x) = i ,  D^2_yU(x,y_j(x)) \mbox{ is positive definite for } j \leq L(x) \},$$ $\mbox{ and }$ $$ F = \mathop{\cup}\limits_{i \geq 1} D_i^\circ$$ with $D_i^\circ$ being the interior of $D_i$. Assume $F^c$ has Lebesgue measure $0$.}
\end{enumerate}

The above assumption is inspired in part by results in parametric
nonlinear programming \cite{JW90}.  It ensures that a modification of Laplace's
method as done in \cite[Theorem 2.1]{H80} applies. We can use it to show that the probability assigned by
$\nu_t^{0, {X}_t}$ to each global minimum $y_i(X_t)$ is proportional
to $(\mbox{Det}\left[ D^2_yU(X_t,y_i(X_t))\right])^{-\frac{1}{2}}$
whenever $X_t \in F$. Though helpful in characterizing the measure it
will still not provide the required regularity to consider $X_t$ as a
classical solution to the differential equation. However, we will be
able to conclude that $X_t$ is a generalized solution to the
differential equation. Towards this we recall a well known concept of
a solution to a differential equation driven by a measurable
function, namely {\em the Filippov solution}.

\begin{definition} Consider the differential equation  given by
\begin{equation}
\frac{d}{dt}{x}(t) = h(x(t)), \ t \geq 0, \ x(0) = x_0, \label{Filip}
\end{equation}
where $h: \mathbb{R}^d \to \mathbb{R}^d$ is a measurable function with at most linear growth. Define the `\textit{enlargement}'
$h_E(\cdot)$ of $h(\cdot)$  to be the set-valued map
\begin{equation}
\label{eqn:hE}
h_E(x) := \mathop{\cap}\limits_{N\subset \R^d: \textsf{Leb}(N) = 0} \,\, \mathop{\cap}\limits_{\delta > 0}\,\, \overline{\mbox{co}}\left(h((x + \delta {\B}_1)\backslash N)\right),
\end{equation}
where \textsf{Leb} denotes Lebesgue measure and $\overline{\mbox{co}}(\cdot)$
denotes the closed convex hull.
An absolutely continuous function $x: [0,\infty) \rightarrow \mathbb{R}^d$ is a Filippov solution to (\ref{Filip}) if it is a solution to the following differential inclusion
\begin{displaymath}
\frac{d}{dt}x(t) \in h_E(x(t)), \ \forall t \geq 0,
\end{displaymath}
with $x(0) = x_0.$
\end{definition}

We refer the reader to \cite{BOQ09} for motivation and various
equivalent definitions of Filippov solutions. {See \cite{BS96} for a comparison of Carath\'{e}odory solutions and Filippov solutions.  Our next result characterizes all limit points
as Filippov solutions of a differential equation.}

Recall the set $F$ from (U4).

\begin{theorem}\label{maintheorem2}
  Assume {\rm (B1), (U1), (U2), (U3), (U4)} and {\rm (S1)}.Then for any sequence
$\varepsilon_n \downarrow 0$ there is a further subsequence,
$\varepsilon_{n_k} \downarrow 0$, along which
$\{X^{\varepsilon_{n_k}}_t, {0 \leq t \leq T}\}$ converges in law on $C([0,T]; \R^d)$ to $\{X_t, 0 \leq t \leq T\}$ which belongs almost surely to the set of Filippov solutions to
\begin{equation}
\frac{d}{dt}X_t =  h(X_t), \forall t \geq 0,  \label{odemain}
\end{equation}
with $X_0 = x_0$ and $h : {\mathbb R}^d \rightarrow {\mathbb R}^d$ is defined almost everywhere as follows:
\begin{equation}\label{defh} h(x) =
\sum_{i=1}^{L(x)} b(x,y_i(x))\frac{\left(\mbox{Det} \left[D_y^2 U(x,y_i(x)) \right] \right)^{-\frac{1}{2}}}{\sum_{j=1}^{L(x)} \left(\mbox{Det} \left[ D_y^2 U(x,y_j(x)) \right] \right)^{-\frac{1}{2}}},
\end{equation}
for all $x \in {F}.$
\end{theorem}

Under (U4), the set $F^c$ has Lebesgue measure 0, and hence the function $h$ is almost everywhere given by the determinantal formula. In general $h$ will not be continuous in $x$, and we will need to consider Filippov solutions of (\ref{odemain}). In some circumstances, however, we may be able to get a classical solution.

As we will see later in the proofs, the measure $\nu_t^{0, X_t}$ in
Theorem \ref{maintheorem1} may depend on the subsequence and
consequently no uniqueness claim is being made about the measure in
Theorem \ref{maintheorem1}. If $\arg\min U(X_t, \cdot)$ is a
singleton then the measure $\nu_t^{0, X_t}$ must be the Dirac measure
on the minimizer. In most other cases Theorem \ref{maintheorem2}
applies.

{\bf Future Directions:} We conclude this section by mentioning a few
possible extensions and open problems. The case when $\alpha >1 $ and
there are no periodicity assumptions for the coupled difussion in
(\ref{ex})-(\ref{wye}) still remains open. So does the case when there
is so called ``full dependence'', when the coefficients in front of
the respective Brownian motions depend on both the slow and the fast
processes. We did not introduce coefficients in front of the driving
diffusion process primarily because we wanted to illustrate the
possible limits when small noise phenomena are present in both time
scales.  Our approach of using nonlinear filtering to characterise
limit points can be generalized to this setting but the spectral gap
estimates for the fast processes which are not reversible will
not be available.

There is a possibility of weakening the assumptions on $U.$ Assumption (U4) imposes a strict behavior of $U(x, \cdot)$ around its global minima. One can try to handle the case when $D^2(x,y_i(x))$ is singular by applying a generalization of Laplace's method (see \cite{AH10}). Further, from the proof of Theorem \ref{maintheorem1} we will be able to infer that, if the rate of convergence of $\|X_t^{\varepsilon_n} - X_t\|$ as $n \to \infty$ is understood, then we can characterize $\nu_t^{0,X_t}$ without assumption (U4).  However, such a rate seems hard to capture given the two timescales and the interdependence of $X^{\varepsilon_n}_t$ on $Y^{\varepsilon_n}_t.$ Towards this an LDP as in \cite{V00} or fluctuation results as in \cite{S14} when $0 <\alpha < 1$ will have to be understood first. Several constants are assumed to be universal in (U1)and (U2), weakening these should be possible and in some cases even our current proof may hold for a restricted set of $\alpha$.

\subsection{Examples}
\label{examples-refinements}

In this section we explore specific examples of $U$ that will help us
understand the assumptions used in Theorem \ref{maintheorem1} and
Theorem \ref{maintheorem2}.

\subsubsection{Weak Convergence and a Classical Solution}

Under assumption (U4), if $L(x) \equiv L$, if $y_i(x)$ were Lipschitz
in $x$ for $1 \leq i \leq L$, and if $F = \mathbb{R}^d$, then $h$is Lipschitz. The ordinary differential equation (\ref{odemain}) is
then well-posed and has a unique solution. We can then strengthen
Theorem \ref{maintheorem2} to say that the process $X^{\varepsilon}$
converges weakly to $X$. We now present an example to illustrate this.

\begin{example} \label{UasU1}
Assume {\rm (S1)} holds for the function $U_1$ given below. Take $m=d=1$, and let $b$ be any function that satisfies {\rm (B1)}. Consider $U_1 : \R^2 \rightarrow \R$ given by
\[  U_1(x,y) =
  \begin{cases}
  y^4 - 2y^2 \left(\frac{1/2 + x^2}{1+x^2}\right) + 1, & |y| \leq 10 \\
  y^4 - 2y^2 + 1, & |y| \geq 20,
  \end{cases}
\]
and for $10 < |y| < 20$ define:
$$
U_1(x,y) := (1-\varrho(|y|)) \left[y^4 - 2y^2 \left(\frac{1/2+x^2}{1+x^2}\right) + 1\right] + \varrho(|y|) \left[y^4 - 2y^2 + 1\right]
$$
with $\varrho:\mathbb{R}_+ \rightarrow [0,1]$ such that $\varrho(|y|)=0$ for $|y| \leq 10$ and $\varrho(|y|) = 1$ for $|y| \geq 20$, and $\varrho(\cdot)$ is a $C^2$ function with both $\varrho'(|y|)$ and $\varrho''(|y|)$ taking the values 0 at $|y| = 10$ and $20$.

It is easy to see that {\rm (U1)} holds. Further, $\nabla_y U_1(x,y)$ and $D^2_y U_1(x,y)$ are continuous for all $x$ and $y$, \eqref{U21} holds for $|y|> R = 20$, and \eqref{U22} holds for a sufficiently large $M$ with $R=20$. To see that \eqref{U23} holds, for $|y| \geq R = 20$, one verifies that the left-hand side of \eqref{U23} is a sixth degree polynomial in $y$ with leading coefficient being negative. Optimizing over $y$ we get the upper bound to be $Ms^2$ for suitably large $M$. Thus one can choose $\eta = 2$ to make \eqref{U23} hold. Hence {\rm (U2)} also holds.

Choose $\varrho$ suitably so that for each $x$, the critical points $y$ satisfying $\nabla_y U_1(x,y) = 0$ also satisfy $|y| \leq 10$. To find the critical points, we may then equate $\nabla_y U_1(x,y) = 4y\left(y^2 - \frac{1/2+x^2}{1+x^2}\right) = 0$. There are then exactly three such points for each $x$. The global minima of $U_1(x,\cdot)$ are then attained at $$y_{1}(x) =  \sqrt{\frac{1/2+x^2}{1+x^2}}, \quad y_{2}(x) = - \sqrt{\frac{1/2+x^2}{1+x^2}}$$ yielding $L(x) = 2$ for all $x$. The point $y=0$ is a local maximum for all $x$.

The quantity $V^1(x)$, by symmetry, is the action functional for moving from $-\sqrt{\frac{1/2+x^2}{1+x^2}}$ to $\sqrt{\frac{1/2+x^2}{1+x^2}}$ and $V^2(x) = 0$. By considering constant velocity paths, it is easy to verify that  action functional is bounded as a function of $x$ and hence Assumption {\rm (U3)} holds.

Finally, $L(x) \equiv 2$ and $D_y^2U_1(x,y_1(x)) = D_y^2U_1(x,y_2(x)) = 8\cdot \frac{1/2+x^2}{1+x^2} \geq 4$ {for all $x \in \R$.} Thus {\rm (U4)} also holds{with $F^c = \emptyset$}.

Theorem \ref{maintheorem2} then implies $X_t$ is a Filippov solution to (\ref{odemain1}) which for this example reduces to
 \begin{equation} \label{odemaine1}\frac{d}{dt}X_t =
\frac{1}{2}b\left(X_t,\sqrt{\frac{1/2+X_t^2}{1+X_t^2}}\right) +
\frac{1}{2}b\left(X_t,-\sqrt{\frac{1/2+X_t^2}{1+X_t^2}}\right), \quad
X_0 = x_0.\end{equation}
Further, from {\rm (B1)}, we note that the
 driving function above is globally Lipschitz. This implies that every
 limit point $X$ is given by the unique classical solution to the
 differential equation (\ref{odemaine1}). Consequently we have that $\{X^\varepsilon_t,~t\in [0,T]\}$ converges in law to the unique solution to the differential equation (\ref{odemaine1}).
\end{example}

\subsubsection{ Merging and Creation of Global Minima}

{We now discuss two illustrative examples where the number of global
minima $L(x)$ varies with $x$. As $x$ varies, global minima may merge
or new global minima may emerge. We begin with an example where global
minima merge. In such an event $D^2_yU(x,y_i(x))$ could have a vanishing
determinant resulting in a nonempty $F^c$ in assumption (U4).  }

\begin{example} \label{UasU2}
Assume {\rm (S1)} holds for the function $U_2$ below. Take $m=d=1$, and let $b$ be any function that satisfies {\rm (B1)}. Similar to Example \ref{UasU1} consider $U_2 : \R^2 \rightarrow \R$ given by
  \[  U_2(x,y) =
  \begin{cases}
  y^4 - 2y^2 \frac{x^2}{1+x^2} + 1, & |y| \leq 10 \\
  y^4 - 2y^2 + 1, & |y| \geq 20,
  \end{cases}
\]
and for $10 < |y| < 20$ define:
$$
U_2(x,y) := (1-\varrho(|y|)) \left[y^4 - 2y^2 \frac{x^2}{1+x^2} + 1\right] + \varrho(|y|) \left[y^4 - 2y^2 + 1\right]
$$
with $\varrho:\mathbb{R}_+ \rightarrow [0,1]$ such that $\varrho(|y|)=0$ for $|y| \leq 10$ and $\varrho(|y|) = 1$ for $|y| \geq 20$, and $\varrho(\cdot)$ is a $C^2$ function with both $\varrho'(|y|)$ and $\varrho''(|y|)$ taking the values 0 at $|y| = 10$ and $20$.

Again, choose $\varrho$ suitably so that for each $x$, the critical points $y$ satisfying $\nabla_y U_2(x,y) = 0$ also satisfy $|y| \leq 10$, and so we may equate $\nabla_y U_2(x,y) = 4y(y^2 - x^2/(1+x^2)) = 0$. The global minimum is then:
\begin{enumerate}
\item[(a)]  attained at $y_{1}(x) =  \frac{x}{\sqrt{1+x^2}}$, $y_{2}(x) = - \frac{x}{\sqrt{1+x^2}}$ when  $x \neq 0$; and
\item[(b)] attained at $y_1(0)=0$ (which is the unique global minimum)  when $x = 0.$
\end{enumerate}
Thus $L(x) = 2$ when $x \neq 0$ and the global minima $y_1(x)$ and $y_2(x)$ merge as $x \rightarrow 0$ yielding $L(0) = 1$.

Following the arguments in Example \ref{UasU1}, we can conclude that {\rm (B1)}, {\rm (U1)}, {\rm (U2)}, {\rm (U3)} hold. Furthermore, $D^2_yU(x,y_i(x))$ is positive definite for all $x \neq 0$ and singular only at $x=0$, {\rm (U4)} also holds {with $F^c = \{0\}$}. From Theorem \ref{maintheorem2}, we know that all limits points are characterized by Filippov solutions to (\ref{odemain}) with
\begin{equation}
\label{h:example2}
h(x) = \frac{1}{2}b\left(x,\frac{x}{\sqrt{1+x^2}}\right) + \frac{1}{2}b\left(x,\frac{-x}{\sqrt{1+x^2}}\right), \mbox{ for all } x \neq 0.
\end{equation}
From  Theorem \ref{maintheorem1} we know that $X$ solves (\ref{odemain1}).

With $L(0) = 1$, we must also have $\nu_t^{0,X_t} = \delta_0$ whenever $X_t = 0$. So we may define $h(0) = b(0,0)$ and we have from assumption {\rm (B1)} that the $h$ in \eqref{h:example2} with $h(0) = b(0,0)$ is a Lipschitz continuous function. As in Example \ref{UasU1}, we obtain convergence in law to the unique solution to the differential equation $$\frac{d}{dt}X_t = \frac{1}{2}b\left(X_t,\frac{X_t}{\sqrt{1+X_t^2}}\right) + \frac{1}{2}b\left(X_t,\frac{-X_t}{\sqrt{1+X_t^2}}\right), \quad  X_0 = x_0.$$
\end{example}

Recall that in Example \ref{UasU1} has $L(x) = 2$ global minima for
all $x$ (i.e. no creation or merging) and Example \ref{UasU2} has
$L(x) = 2$ global minima for all $x \neq 0$, but they merge as $x
\rightarrow 0$ to give $L(0) = 1$. In the following example, we
consider a different variation where a new global minimum is created.

\begin{example} \label{UasU3}
Consider
$$U_3(x,y) = U_1(x,y) + \phi(x) y^4 1\{ y \geq 0 \},$$
where $U_1$ is as in Example \ref{UasU1} and $\phi(x)$ is any smooth and strictly increasing function that is strictly positive when $x > 0$, equals 0 when $x = 0$, strictly negative when $x < 0$, and $\phi(x) \geq -1/2$ for all $x$. Note that this is a perturbation of $U_1$. When $x > 0$, the perturbation term $\phi(x) y^4 1\{ y \geq 0 \}$ lifts the graph $U(x,\cdot)$ for $y > 0$ but leaves it unchanged for $y \leq 0$, and therefore the left minimum of $U_1(x,\cdot)$ is the unique global minimum of $U_3(x,\cdot)$. Similarly, when $x < 0$, the perturbation pushes the graph gently down for $y > 0$, leaves it unchanged for $y \leq 0$, and therefore the unique global minimum of $U_3(x,\cdot)$ is strictly positive. When $x=0$ however, we get $U_3(0,\cdot) = U_1(0,\cdot)$ and we therefore have two global minima.

Thus $L(x) = 1$ for all $x \neq 0$, $L(0) = 2$, and assumption (U4) holds with $F^c = \{0\}$. It is easy to see that all assumptions for Theorem \ref{maintheorem2} hold and Theorem \ref{maintheorem2} applies. However, if $b(0,1/\sqrt{2}) \neq b(0,-1/\sqrt{2})$, the resulting $h$ in (\ref{defh}) has $h(0-) \neq h(0+)$. So we will not in general have a classical solution to (\ref{odemain}), but we do have a generalized solution, namely, the Filippov solution.

More generally, for any nonconstant $b(\cdot,\cdot)$, one can choose $0 \leq \phi(\cdot) \leq 1/2$ arising from a Lipschitz-continuous distance function with distance taken from a suitable generalized Cantor-type set, so Theorem \ref{maintheorem2} applies. In this case as well the nature of $h$ will be such that we can at best ensure that all limit points are generalized Filippov solutions to (\ref{odemain}).

{An example of $U(\cdot, \cdot)$ that does not satisfy (U4) is
$U_4(x,y) = \phi(x) y^2 + y^4,$
where $\phi^{-1}(0)$ has positive Lebesgue measure.}
\end{example}

{\bf Layout of the paper:} The rest of the paper is organized as follows. In Section \ref{pomt} we present three key results: Proposition \ref{p:spgp} (establishes a spectral gap bound for the rescaled fast process (\ref{associatedsde})), Proposition \ref{p:tightness} (identifies limit points), and Proposition \ref{characterisation} (characterizes a given limit point). These are used in the proof of the main results. Section \ref{pop:spgp} is devoted to the proof of Proposition \ref{p:spgp}, Section \ref{pop:tightness} is devoted to the proof of Proposition \ref{p:tightness}, and Section \ref{pop:characterisation} is devoted to the proof of Proposition \ref{characterisation}. In the appendix, in Appendix \ref{existenceanduniqueness} we prove some
auxiliary results concerning existence and uniqueness of the slow-fast small noise diffusions (\ref{ex}) and (\ref{wye}). In Appendix \ref{app:nonlinearfiltering}, we provide a general result in Proposition \ref{fkksde} on the Fujisaki-Kallianpur-Kunita nonlinear filtering equation satisfied by the conditional law of one component of a pair of coupled diffusions given the other. We conclude the article with an extension of Laplace's principle in Appendix \ref{appendix:Laplace}.

{\bf Convention on constants:}  We shall now
fix $\alpha$ and the functions $b,U,s$. Unless otherwise mentioned, we
assume that (B1), (U1), (U2), (U3) and (S1) are satisfied as stated above
with the associated constants. All other positive valued constants
whose values are not important will be denoted by $c_1, c_2,\ldots,$
and their dependencies on parameters if needed will be mentioned
inside parentheses, e.g.,  $c_1(\alpha).$ For such constants, the numbering
will begin afresh in each new result and proof.

\section{Proof of Theorem \ref{maintheorem1} and Theorem \ref{maintheorem2}} \label{pomt}

In this section we shall state three key propositions and prove Theorem \ref{maintheorem1} and Theorem \ref{maintheorem2}. The proofs of the propositions follow in subsequent sections.

{Our approach is inspired by the foundations laid in \cite{FW12}.
  The fast moving $Y_t^{\varepsilon}$ process approaches its
  stationary distribution $\nu^{\varepsilon, X_t^{\varepsilon}}$ and
  this in turn approaches a limiting measure ``$\nu_t^{0,X_t}$'' as $\varepsilon \rightarrow 0$. Thus the slow
  process $X^{\varepsilon}_t$ as $\varepsilon \rightarrow 0$ will now
  observe an averaging principle in $Y^{\varepsilon}_t$ as determined
  by $\nu_t^{0,X_t}$. To make the above rigorous we will need to
  quantify to what extent $Y^{\varepsilon}_t$ has equilibrated to
  ``$\nu^{\varepsilon, X_t^{\varepsilon}}$'' along with the rate of
  convergence of ``$\nu^{\varepsilon, X_t^{\varepsilon}}$'' to
  $\nu_t^{0,X_t}$ as $\varepsilon \rightarrow 0$.  However implementing this program of analysis turns out to be delicate
  due to the presence of small noise limit dictated by
  $s(\varepsilon)$. We will see this manifest itself in the spectral  gap estimate for the fast process, which we will establish first.}

Fix $t \in [0,T], x \in \R^d, y \in \R^m, s \geq t$. Consider the
stochastic differential equation
\begin{equation} \label{associatedsde}
  Z^{t,\varepsilon,x}_s \ = \ y - \int^s_t \nabla_y U(x, Z^{t,\varepsilon, x}_r) dr +  s(\varepsilon) (\overline{W}_s - \overline{W}_t), s \geq t,
\end{equation}
with $\overline{W}_t$ being a Brownian motion. One may view the above
stochastic differential equation as being obtained from (\ref{wye}) by
first freezing $X^{\varepsilon} \equiv x$, then scaling time by
$\varepsilon$ and setting $\overline{W}_t :=
\frac{1}{\sqrt{\varepsilon}} W_{\varepsilon t}, t\geq 0.$ The small
noise limit in (\ref{associatedsde}) (i.e. $s(\varepsilon) \rightarrow
0$ as $\varepsilon \rightarrow 0$) has been well studied in the
literature. Hwang and Sheu \cite{HS90} gave explicit decay rates
for the second eigenvalue of the Fokker-Planck operator associated
with the generator of (\ref{associatedsde}) and provided connections to
simulated annealing (where the exact formulation of $s(\varepsilon)$ can
be identified). The small noise phenomenon in (\ref{associatedsde})
can be used to identify the global minima of the function $U$ and has applications in simulated annealing (see \cite{HS90}). These
and the other physical phenomenon of metastability have been explored by
Bovier et al. in \cite{ BEGK00, BEGK01, BEGK04} and by Eckhoff in
\cite{E05}. Recently in \cite{BB09}, limits of invariant measures of
(\ref{associatedsde}) under the small noise limit were understood via a
control theoretic approach.

For $ s \geq t, f \in C^2_b
(\mathbb{R}^m)$, consider the Feller semigroup of
the process $Z^{t,\varepsilon,x}_{s}$ defined by
\[
T^{t,\varepsilon, x}_s f(y) \ = \ \E_y \Big[ f(Z^{t,\varepsilon,x}_{s})\Big] \
\]
with the corresponding generator given by
\begin{equation}\label{generator}
\LA^{\varepsilon,x} f (y) = \frac{s(\varepsilon)^{2}}{2}\Delta f(y) - \langle \nabla_yU(x, y),
\nabla f(y) \rangle.
\end{equation}

Our first proposition describes the invariant measure of $Z^{t, \varepsilon, x}$ and provides a uniform rate of convergence to stationarity using a spectral gap estimate.

\begin{proposition} \label{p:spgp} {\bf (Spectral Gap Estimate)}
 Let $x \in \R^d, y \in \R^m$ and $t \in [0,T]$.
  \begin{enumerate}
    \item[(a)] The stochastic differential equation (\ref{associatedsde}) has a unique strong solution  equipped with a unique invariant probability measure $\nu^{\varepsilon,x}(dy)$ given by
 \begin{equation}\label{epsiloninvariantmeasure}
 \nu^{\varepsilon,x}(dy) :=  C(\varepsilon, x)^{-1} e^{-\frac{2U(x, y)}{s(\varepsilon)^2}}dy,
 \end{equation}
 where $ 0 < C(\varepsilon, x) < \infty$ is the normalizing factor.

 \item[(b)] Fix $\delta > 0$. For all sufficiently small $\varepsilon$, there exists a $c_1 > 0$ such that for all $s > { t + 1}$ and
$f \in C^2_b (\mathbb{R}^m)$
\begin{equation} \label{spinf}
\parallel T^{t,\varepsilon,x}_s f - \nu^{\varepsilon,x} (f) \parallel_\infty \quad \leq  \quad \parallel f\parallel_{\infty} \,
{ e^{\frac{c_1}{s(\varepsilon)^2} -(s-t) \exp\left({-\frac{(\Lambda+\delta)}{s(\varepsilon)^2}}\right)}},
\end{equation}
where $0 \leq \Lambda < \infty$ is as in (\ref{u3}).
\end{enumerate}
\end{proposition}

\begin{remark}
\label{rem:mixing}
The spectral gap for reversible diffusion (\ref{associatedsde}) is proved in \cite[Theorem 3.1]{HS90} and from this (\ref{spinf}) will follow in the $L_2$ sense. We however need the estimate in the infinity norm and the spectral gap to be independent of $x \in \R^d.$ These are achieved respectively by ultracontractivity due to (\ref{U23}) of  assumption (U2) and (\ref{u3}) of assumption (U3) resulting in an extra factor $\exp\left\{ \frac{c_1}{s(\varepsilon)^2} \right\}$.

As we see later, we will choose $s-t$ to be $\varepsilon^{-\theta}$ for some $\theta > 0$. Hence $s(\varepsilon)$ as in \eqref{svarepsilon} of Assumption (S1) ensures that the process has mixed and the right-hand side of \eqref{spinf} goes to zero.
\end{remark}

{}

Our next step is to establish tightness of $X^\varepsilon$ along with tightness of conditional laws of $Y^\varepsilon $ given  $X^\varepsilon$ in a specific topology. For $s > 0$, set ${\mathcal F}^\varepsilon_s$ as the
$\P$-completion of $\cap_{s' > s} \sigma(X^\varepsilon_{u} , u \leq
s').$ Define $\pi^\varepsilon_s \in \PA(\mathbb{R}^m)$ via
\begin{equation}\label{conditionalLaw}
\pi^\varepsilon_s(f) := \ \E [f(Y^\varepsilon_s) |
 {\mathcal F}^\varepsilon_s], \quad \forall f \in C_b(\R^m).
\end{equation}
Using \cite[Theorem 4.1]{Wo71}, one can rewrite (\ref{ex}) in the form
\begin{equation}\label{ex1}
  X^\varepsilon_t \ = x_0 + \int_0^t \int b(X^\varepsilon_s, y) \pi^\varepsilon_s(dy) ds + \varepsilon^\alpha \eta^{\varepsilon}_t,
\end{equation}
where $\eta^{\varepsilon}_t$
is an $\mathbb{R}^d$-valued Wiener process under $\P$.
 Let $\overline{\mathbb{R}}^m$ denote the one point compactification of
$\mathbb{R}^m$. We equip
$$
\widetilde{\PA} = \{ \zeta : ~(\zeta:[0,
  T] \rightarrow \PA (\overline{\mathbb{R}}^m)) \mbox{ and is
  measurable} \}
$$
with the coarsest topology that renders continuous
the maps $$\zeta \in \widetilde{\PA} \rightarrow \int^s_u g(a) \int
f(y) \zeta_a(dy) da$$ for all $ 0\leq u < s \leq T$, $g \in L^2 [u,
  s], f \in C(\overline{\mathbb{R}}^m)$.  We will view $(X^{\varepsilon}, \pi^{\varepsilon}) := (X^{\varepsilon}_t, \pi^{\varepsilon}_t)_{t \in [0,T]}$ as elements of $C([0, T] ; \mathbb{R}^d) \times \widetilde{\PA}$.

{The above approach towards topologizing the path space of the
  conditional density of $Y^{\varepsilon}$ is borrowed from the
  relaxed control framework in control theory. This is described in
  Chapter 2 of \cite{ABG12}. More specifically, the topology is compact and metrizable as explained in \cite[Section 2.3]{ABG12}.   Our next proposition asserts
  tightness and identifies a limit point with which we will work.}

\begin{proposition} \label{p:tightness} {\bf (A limit point)}
The laws of $\{(X^\varepsilon, \pi^\varepsilon):  0 < \varepsilon < 1\}$ are tight in the space $\PA ( C(
[0, T] ; \mathbb{R}^d) \times \widetilde{\PA})$. Further, there exists a sequence
  $\varepsilon_n  \rightarrow 0$ as $n \rightarrow \infty$ such that
\begin{enumerate}
\item[(a)]  $({X}^{\varepsilon_n}, {\pi}^{\varepsilon_n}) \rightarrow ({X}, {\pi})$ weakly as $n \rightarrow \infty,$

\item[(b)] there exists a filtered probability space $(\tilde{\Omega}, \tilde{\FA}, \tilde{\P})$,  random processes \\$(\tilde{X}^{\varepsilon_n},\tilde{\pi}^{\varepsilon_n}, \tilde{\eta}^{\varepsilon_n}) \stackrel{d}{=} (X^{\varepsilon_n}, \pi^{\varepsilon_n}, \eta^{\varepsilon_n})$  and $(\tilde{X},\tilde{\pi},\tilde{\eta}) \stackrel{d}{=} (X, \pi, \eta)$ such that
       \ben
       \label{eqn:second-moment-convergence}
(\tilde{X}^{\varepsilon_n}, \tilde{\pi}^{\varepsilon_n}, \tilde{\eta}^{\varepsilon_n}) \rightarrow ( \tilde{X}, \tilde{\pi}, \tilde{\eta}) \mbox{ a.s.,  and } \,\,
      \tilde{\E} \left[ \sup_{t \in [0,T]}\norm{ \tilde{X}^{\varepsilon_n}_t -  \tilde{X}_t }^2 \right] \rightarrow 0
       \een
   as $n \rightarrow \infty$.
   \end{enumerate}
\end{proposition}

{From the above result we have a candidate limit point for
  $X^{\varepsilon}$ and a limit point for the conditional density
  $\pi^\varepsilon$ in the $\epsilon \rightarrow 0$ limit. As
  indicated earlier we will use filtering theory to understand the
  limit point of $\pi^\varepsilon$. The chosen topology enables the
  use of the spectral gap estimate to identify how $s(\varepsilon)$ should decay to $0$ as $\varepsilon \rightarrow 0$ in order to
  establish that any limit point of $\pi^\varepsilon$ coincides with a probability measure supported on the $\arg \min \{U(X_t,\cdot)\}$.}

One could directly show tightness of $Y^\varepsilon$ but
characterizing the limit point does not seem to be straightforward
(except in the case when $U(x, \cdot)$ has a unique global
minimum).{ However, the above leads to a much simpler approach to
  the averaging result because it enables us to avoid reliance on empirical measures of the fast process (which are more difficult to handle). Further, as discussed in the introduction, the
  probability measure-valued process of conditional laws has its own
  well defined evolution given by the Fujisaki-Kunita-Kallianpur
  equation of nonlinear filtering ( see Proposition \ref{fkksde}). This facilitates the characterization of its weak limit points in a straightforward manner, which is our next result.}

\begin{proposition} {\bf (Characterization of $\pi$)} \label{characterisation}
  Let $\gamma = \min\{1-\alpha, \frac{1}{2}\}$. Let $\varepsilon_n >0, \tilde{X}, \tilde{\pi}$ be as constructed in Proposition \ref{p:tightness}.  There exists a subsequence $\varepsilon_{n_k}$  such that for all $f \in C_0^2(\R^m)$:
  \begin{enumerate}
    \item [(a)]      \ben  \label{l3a}\lim_{k \rightarrow \infty} \frac{1}{\varepsilon_{n_k}^{\gamma}} \int^{t + \varepsilon_{n_k}^{\gamma}}_t\tilde{\pi}^{\varepsilon_{n_k}}_s (f) ds  = \tilde{\pi}_t(f)
      \een
      $\mbox{ for almost every } t \in [0,T], \mbox{almost surely}$;
    \item[(b)] \ben \label{l3} \lim_{k \rightarrow \infty} \left| \nu^{\varepsilon_{n_k},\tilde{X}^{\varepsilon_{n_k}}_t}(f)-  \frac{1}{\varepsilon_{n_k}^{\gamma}} \int^{t + \varepsilon_{n_k}^{\gamma}}_t{\tilde{\pi}}^{\varepsilon_{n_k}}_s (f) ds   \right| = 0,
      \een
      for  all $t \in [0,T]$ almost surely; and
    \item[(c)] Almost surely, for almost every $t \in [0,T]$, $\nu^{\varepsilon_{n_k}, \tilde{X}^{\varepsilon_{n_k}}_t}$ converges weakly to a probability measure $\nu_t^{0, \tilde{X}_t}$ supported on $\arg\min U(\tilde{X}_t, \cdot)$ and further $\tilde{\pi}_t  = \nu_t^{0,\tilde{X}_t}$.
    \item[(d)] If (U4) holds and if $\tilde{X}_t \in F$, then the measure $\nu_t^{0, \tilde{X}_t}$ from (c) is given by
        $$ \nu_t^{0, \tilde{X}_t}(\cdot) =  \sum_{i=1}^{L(\tilde{X}_t)} b(x,y_i(\tilde{X}_t))\frac{\left(\mbox{Det}\left [D_y^2 U(x,y_i(\tilde{X}_t))\right]\right)^{-\frac{1}{2}}}{\sum_{j=1}^{L(\tilde{X}_t)} \left(\mbox{Det}\left [D_y^2 U(x,y_j(\tilde{X}_t))\right] \right)^{-\frac{1}{2}}}.$$
      \end{enumerate}

\end{proposition}

The above proposition contains the main architecture of the proof of
Theorem \ref{maintheorem1}. It works with the sequence
$\{ \varepsilon_n\}$ and the associated limit point from Proposition
\ref{p:tightness}. Part (a) shows that convergence of the conditional
densities holds in the small time-averaged limit along a
subsequence. The topology borrowed from \cite{ABG12} is made use of in
this step.  Part (b) contains the key step that is used to understand
the two ``limits'', first one in which the fast process approaches
stationarity resulting in the averaging phenomenon and the second one
in which the stationary measure approaches its limit due to the
presence of small noise in (\ref{wye}). In Proposition
  \ref{l:ke} we show a second moment estimate. It is here that we critically
  benefit from the filtering theory approach, understand the role played by the decay rate of $s(\varepsilon)$ to zero as
  $\varepsilon \rightarrow 0$, and observe the need to
  choose $C$ large enough to achieve the result.

  Part (c) characterizes
  all subsequential weak limits of $\nu^{\varepsilon,
    \tilde{X}^{\varepsilon}_t}$ as measures supported on $\arg\min
  U(\tilde{X}_t,\cdot)$ denoted by $\nu_t^{0, \tilde{X}_t}$. This confirms that the measures $\tilde{\pi}_t$, known to be supported on $\overline{\mathbb{R}}^m$, are actually supported on $\arg\min
  U(\tilde{X}_t,\cdot)$ for almost every $t \in [0,T]$. Finally,
  in Part (d), assumption (U4) is used to enable the implementation of Laplace's principle to arrive at a determinantal formula for
  subsequential limits.

  { We note that the characterization of $\tilde{\pi}$ may change with
    the choice of the subsequence taken in the previous parts and
    consequently there is no uniqueness claim being made about the
    { measure $\nu_t^{0, \tilde{X}_t},$ under (U1)-(U3) alone. Of
      course, if $\arg\min U(\tilde{X}_t, \cdot)$ is a singleton, it
      is perforce unique, being the Dirac measure on the minimizer. In
      part this motivated assumption (U4) under which a modification
      of the Laplace's method holds and the probability assigned by
      $\nu_t^{0, \tilde{X}_t}$ to each global minima is      proportional to $\left(\mbox{Det}\left[
        D^2U(\tilde{X}_t,y_i(\tilde{X}_t))\right]\right)^{-\frac{1}{2}}$ provided the Hessian (in $y$) of $U$ at all global minima of
      $U(\tilde{X}_t, \cdot)$ are {\em uniformly positive definite}  in a neighborhood of $\tilde{X}_t$. }}

We are now ready to present the proofs of Theorem \ref{maintheorem1}
and Theorem \ref{maintheorem2}. We will begin by setting up common
notation required for both and will then present the proof of each. From (\ref{ex1}) we have that
\begin{eqnarray}\label{eq1theorem2.1}
X^{\varepsilon}_t & = &  x_0 + \int^t_0 \int b(X^\varepsilon_s, y) \pi^\varepsilon_s (dy) ds
+ \varepsilon^{\alpha} \eta^\varepsilon_t, \quad t \in [0,T].\nonumber
\end{eqnarray}
Let $\varepsilon_n \rightarrow 0$ denote the subsequence identified in Proposition \ref{characterisation}. So there exist a probability space and processes $(\tilde{X}^{\varepsilon_n}, \tilde{\pi}^{\varepsilon_n}, \tilde{\eta}^{\varepsilon_n}, \tilde{X}, \tilde{\pi}, \tilde{\eta})$ such that
\begin{itemize}
  \item $(X^{\varepsilon_n}, \pi^{\varepsilon_n}, \eta^{\varepsilon_n})$ and $(\tilde{X}^{\varepsilon_n}, \tilde{\pi}^{\varepsilon_n}, \tilde{\eta}^{\varepsilon_n})$ have the same law for $n \geq 1$;
  \item $(X, \pi, {\eta})$ and $(\tilde{X}, \tilde{\pi}, \tilde{\eta})$ have the same law;
  \item $\tilde{X}^{\varepsilon_n} \rightarrow \tilde{X}$ and $\tilde{\eta}^{\varepsilon_n} \rightarrow \tilde{\eta}$ in $C([0,T];\R^d)$, and $\tilde{\pi}^{\varepsilon_n} \rightarrow \tilde{\pi}$ in $\tilde{\PA}$, a.s.
\end{itemize}

Set $\tilde{\xi}^{\varepsilon_n}_t := \tilde{X}^{\varepsilon_n}_t - \varepsilon_n^{\alpha} \tilde{\eta}^{\varepsilon_n}_t$. Then,
\begin{eqnarray} \label{comthmstep1}
  \tilde{\xi}^{\varepsilon_n}_t & = & \tilde{X}^{\varepsilon_n}_t - \varepsilon_n^{\alpha} \tilde{\eta}^{\varepsilon_n}_t
   = x_0 + \int^t_0 \int b(\tilde{X}^{\varepsilon_n}_s, y) \tilde{\pi}^{\varepsilon_n}_s (dy) ds,
\end{eqnarray}
Since $\tilde{X}^{\varepsilon_n} \rightarrow \tilde{X}$ and $\tilde{\eta}^{\varepsilon_n} \rightarrow \tilde{\eta}$ in $C([0,T];\R^d)$, a.s., we have
\begin{equation}
 \label{comthmstep2} \tilde{\xi}^{\varepsilon_n} \rightarrow \tilde{X}  \mbox{ in $C([0,T]; \R^d)$, a.s.}
 \end{equation}
{For all $s \in [0,T]$, define
\begin{equation} \label{dens}
 \delta_{n,s}:= \left\| \int b(\tilde{X}^{\varepsilon_n}_s,y) \tilde{\pi}^{\varepsilon_n}_s(dy) - \int b(\tilde{X}_s,y) \tilde{\pi}^{\varepsilon_n}_s(dy) \right\|
 \end{equation}
and define
\begin{equation} \label{tens}
  \tau_{n,s} := \left \| \int b(\tilde{X}_s,y) \tilde{\pi}^{\varepsilon_n}_s (dy) - \int b(\tilde{X}_s,y) \tilde{\pi}_s(dy) \right\|.
\end{equation}
By Proposition \ref{characterisation}(c), $\tilde{\pi}_s = \nu_s^{0,\tilde{X}_s}$ for almost every $s \in [0,T]$, a.s., with $\nu_s^{0,\tilde{X}_s}$ being  a probability measure supported on $\arg\min U(\tilde{X}_s, \cdot)$.  We can therefore write
\begin{equation} \label{tens1}
  \tau_{n,s} = \left \| \int b(\tilde{X}_s,y) \tilde{\pi}^{\varepsilon_n}_s (dy) - \int b(\tilde{X}_s,y) \nu_s^{0,\tilde{X}_s}(dy) \right\| \mbox{ for almost every } s \in [0,T].
\end{equation}}

{\em Proof of Theorem \ref{maintheorem1}}: {Observe that by (\ref{comthmstep1}), (\ref{dens}), and (\ref{tens1}),  with some simple algebra we have}
\begin{equation}\label{mt1step1}
\| \tilde{X}_t - x_0 -  \int^t_0 \int b(\tilde{X}_s, y) \nu_s^{0,\tilde{X}_s}(dy)ds\| \leq  \| \tilde{X}_t -   \tilde{\xi}^{\varepsilon_n}_t\| +\int_0^t \delta_{n,s} ds + \int_0^t \tau_{n,s}ds.
\end{equation}
By the Lipschitz property of $b$ in (B1), we  have for all $s \in [0,T]$,
\[
\int_0^t \delta_{n,s} ds\leq K \int_0^t\|\tilde{X}^{\varepsilon_n}_s - \tilde{X}_s\|ds \leq K T \sup_{s \in [0,T]} \|\tilde{X}^{\varepsilon_n}_s - \tilde{X}_s\|
 \]
 and so  \begin{equation}\label{mt1step2}
   \int_0^t \delta_{n,s}ds  \rightarrow 0 \mbox{ for all $t \in [0,T]$ a.s.}
   \end{equation}
 By Proposition \ref{p:tightness}(b), as noted earlier,
 $\tilde{\pi}^{\varepsilon_n} \rightarrow \tilde{\pi}$ in
 $\tilde{\PA}$. By Proposition \ref {characterisation}(c), $\tilde{\pi}_s$ is supported on $\arg\min\{U(\tilde{X}_s,\cdot)\}$ which is a finite set for each $s \geq 0$. Consequently, using the topology on $\tilde{\PA}$ it is standard to see that for $h \in C_b (\R^m)$
      $$   \left \| \int h(y) \tilde{\pi}^{\varepsilon_n}_s (dy) - \int  h(y) \tilde{\pi}_s(dy) \right\| \rightarrow 0 \mbox{ almost every  } s \in [0,T].$$
As $b(\tilde{X_s}, \cdot)$ is a bounded (though random) continuous function, we then have
     $$   \left \| \int b(\tilde{X}_s,y) \tilde{\pi}^{\varepsilon_n}_s (dy) - \int  b(\tilde{X}_s,y) \tilde{\pi}_s(dy) \right\| \rightarrow 0 \mbox{ almost every  } s \in [0,T].$$

By (\ref{tens}) and (\ref{tens1}), this is the same as
$\tau_{n,s} \rightarrow    0$ for almost every $s \in [0,T]$, a.s. An application of the
     dominated convergence theorem then yields that \begin{equation} \label{mt1step3} \int_0^t \tau_{n,s}ds \rightarrow 0 \mbox{ for all } t \in [0,T] \mbox{ a.s.}
     \end{equation}
     So using (\ref{mt1step1}), and by (\ref{comthmstep2}), (\ref{mt1step2}), and (\ref{mt1step3}) we have
     \[\tilde{X}_t = x_0 +  \int^t_0 \int b(\tilde{X}_s, y)  \nu_s^{0,\tilde{X}_s}(dy)ds\]
     for all $ t \in [0,T]$  a.s.
   {This completes the proof.}
\qed

The method of proof for Theorem \ref{maintheorem2} is adapted from
Theorem 4 in \cite{BOQ09} with some key differences. We present it next.

{\em Proof of Theorem \ref{maintheorem2}}:  From Proposition \ref{p:tightness}(b), we have
\begin{equation}
  \label{eqn:L2-convergence}
  \tilde{\E}\left[ \sup_{t \in [0,T]} \| \tilde{X}^{\varepsilon_n}_t - \tilde{X}_t \|^2 \right] \rightarrow 0 \quad \mbox{ as } n \rightarrow \infty.
\end{equation}
Since
\[
\| \tilde{\xi}^{\varepsilon_n}_t - \tilde{X}_t \|^2 = \| \tilde{X}^{\varepsilon_n}_t - \varepsilon_n^{\alpha} \tilde{\eta}^{\varepsilon_n}_t - \tilde{X}_t \|^2 \leq 2 \| \tilde{X}^{\varepsilon_n}_t  - \tilde{X}_t \|^2 + 2 \varepsilon_n^{2\alpha} \| \tilde{\eta}^{\varepsilon_n}_t \|^2,
\]
this together with the facts
$\tilde{\E}[\sup_{t \in[0,T]} \| \tilde{\eta}^{\varepsilon_n}_t \|^2] < \infty$, $\varepsilon_n \rightarrow 0$, and (\ref{eqn:L2-convergence}) yields
\begin{equation}
  \label{eqn:y-tilde-bound}
  \tilde{\E}\left[ \sup_{t \in [0,T]} \| \tilde{\xi}^{\varepsilon_n}_t - \tilde{X}_t \|^2 \right] \rightarrow 0 \quad \mbox{ as } n \rightarrow \infty.
\end{equation}
Observe now that since
\begin{eqnarray*}
  \tilde{\xi}^{\varepsilon_n}_t & = & \tilde{X}^{\varepsilon_n}_t - \varepsilon_n^{\alpha} \tilde{\eta}^{\varepsilon_n}_t
   = x_0 + \int^t_0 \int b(\tilde{X}^{\varepsilon_n}_s, y) \tilde{\pi}^{\varepsilon_n}_s (dy) ds,
\end{eqnarray*}
we can write
\[
  \frac{d}{dt}\tilde{\xi}^{\varepsilon_n}_t = \int b(\tilde{X}^{\varepsilon_n}_t, y) \tilde{\pi}^{\varepsilon_n}_t(dy),
\]
and in view of the boundedness of $b$ in Assumption (B1), there is a $0< c_1 < \infty$ such that
\begin{eqnarray}
  \label{eqn:y-tilde-derivative-bound}
  \tilde{\E}\left[ \sup_{t \in [0,T]} \left\| \tilde{\xi}^{\varepsilon_n}_t \right\|^2 \right] \leq c_1 \quad \mbox{ and } \quad
  \tilde{\E}\left[ \sup_{t \in [0,T]} \Big\| \frac{d}{dt} \tilde{\xi}^{\varepsilon_n}_t \Big\|^2 \right] \leq c_1, \quad \forall n \geq 1.
\end{eqnarray}
In view of (\ref{eqn:y-tilde-bound}) and (\ref{eqn:y-tilde-derivative-bound}), there is a subsequence that converges weakly in the space
\[
  W^{1,2} := \left\{ Z \in L^2([0,T] \times \tilde{\Omega}; \R^d), Z' \in L^2([0,T] \times \tilde{\Omega}; \R^d) \right\},
\]
that is, there is some process $U$ such that
\begin{eqnarray*}
  \tilde{\xi}^{\varepsilon_n} & \rightarrow & \tilde{X} \quad \mbox{ in } L^2, \\
  \tilde{\E} \left[ \int_0^T \frac{d}{dt}\tilde{\xi}^{\varepsilon_n}_t \phi(t) dt \right] & \rightarrow & \tilde{\E} \left[ \int_0^T U_t \phi(t) dt \right],
\end{eqnarray*}
for any (nonrandom) $\phi \in W^{1,2}$. We next argue that $U = \frac{d}{dt}\tilde{X}$. Integrating the left-hand side above by parts, we get
\begin{eqnarray*}
&\tilde{\E} \left[ \tilde{\xi}^{\varepsilon_n}_T \phi(T) - x_0 \phi(0) \right] - \tilde{\E} \left[ \int_0^T \tilde{\xi}^{\varepsilon_n}_t ~\frac{d}{dt}\phi(t) dt \right]  ~\rightarrow~ \tilde{\E} \left[ \tilde{X}_T \phi(T) - x_0 \phi(0) \right] - \tilde{\E} \left[ \int_0^T \tilde{X}_t ~\frac{d}{dt}\phi(t) dt \right],
\end{eqnarray*}
whence
\[
\tilde{\E} \left[ \int_0^T \frac{d}{dt}\tilde{X}_t \phi(t) dt \right] = \tilde{\E} \left[ \int_0^T U_t \phi(t)dt \right].
\]
Since $\phi \in W^{1,2}$ was arbitrary, with the only restriction that it is nonrandom, we have established that $U_t = ~\frac{d}{dt}\tilde{X}_t$ for almost every $t \in [0,T]$, a.s. Thus $\tilde{\xi}^{\varepsilon_n} \rightarrow  \tilde{X}$ weakly in $W^{1,2}$.

Recall definition of $\delta_{n,t}$ and $\tau_{n,t}$ from (\ref{dens}) and (\ref{tens}), respectively. Let $\delta_n\ := \sup_{t \in [0,T]} \|\tilde{X}^{\varepsilon_n}_t - \tilde{X}_t\|$. As discussed earlier, by the Lipschitz property of $b$ in (B1), we then have for all $t \in [0,T]$,
\[
 \delta_{n,t} \leq T K \delta_n
\]

{Using (\ref{dens}), (\ref{tens}), (\ref{tens1}) and  the triangle inequality we see that the derivative $\frac{d}{dt}\tilde{\xi}^{\varepsilon_n}$ satisfies, for almost every $t \in [0,T]$,}
\begin{eqnarray}
  \frac{d}{dt}\tilde{\xi}^{\varepsilon_n}_t &=&  \int b(\tilde{X}^{\varepsilon_n}_t,y) \tilde{\pi}_t^{\varepsilon_n}(dy)  \nonumber \\
  &\in & \int b(\tilde{X}_t,y) \nu_t^{0,\tilde{X}_t}(dy) + (\tau_{n,t} + \delta_{n,t}) \bar{\mathbb{B}}_1  \nonumber \\ &\subset & \int b(\tilde{X}_t,y) \nu_t^{0,\tilde{X}_t}(dy)  + (\tau_{n,t} + K \delta_n) \bar{\mathbb{B}}_1 \nonumber \\
  & = & \int b(\tilde{X}_t,y) \nu_t^{0,\tilde{X}_t}(dy) + \gamma_{n,t} \bar{\mathbb{B}}_1 \quad \mbox{(where $\gamma_{n,t} = \tau_{n,t} + K\delta_n$)}.
  \label{mt2s1}
\end{eqnarray}
Let  $h(\cdot)$ be as  defined in (\ref{defh}). By assumption (U4) and  Proposition \ref{characterisation}(d), $\mbox{ whenever } \tilde{X}_t \in F,$  we have that $$ \int b(\tilde{X}_t,y) \nu_t^{0,\tilde{X}_t}(dy)= \sum_{i=1}^{L(\tilde{X}_t)} b(\tilde{X}_t,y_i(\tilde{X}_t))\frac{\left(\mbox{Det}\left [D_y^2 U(\tilde{X}_t,y_i(\tilde{X}_t))\right]\right)^{-\frac{1}{2}}}{\sum_{j=1}^{L(\tilde{X}_t)} \left(\mbox{Det}\left [D_y^2 U(\tilde{X}_t,y_j(\tilde{X}_t))\right] \right)^{-\frac{1}{2}}} = h(\tilde{X}_t). $$
Now consider the enlargement $h_E$ of $h$ defined in (\ref{eqn:hE}) as the smallest upper semi-continuous set-valued map with closed convex values such that $h(x) \in h_E(x)$ for almost all $x \in \R^d$.

Define $f,g: \R_+ \times \R^d \rightarrow \R$ by
$$f(t,x) = h(x) \mbox{  and } g(t,x) = \left \{\begin{array}{ll}\int b(\tilde{X}_t,y) \nu_t^{0,\tilde{X}_t}(dy)& \mbox{ if } \tilde{X}_t =x \mbox{ and } x \in F\\
h(x)  & \mbox{otherwise} \end{array} \right. $$
for all $(t,x) \in \R_+ \times \R^d.$ We know that $g =f $ a.e on $\R_+ \times \R^m$ and consequently by \cite[Proposition 2(ii)]{BOQ09} we have  $g_E =f_E.$ As $f$ does not depend on $t$ it is easy to see that the enlargement  $f_E(t,x) = h_E(x)$ for all $t \in [0,T]$ and $x \in \mathbb{R}^d$. Therefore, from (\ref{mt2s1}), we have
\begin{eqnarray*}
  \frac{d}{dt}\tilde{\xi}^{\varepsilon_n}_t &\in & h_E(\tilde{X}_t) + \gamma_{n,t} \bar{\mathbb{B}}_1, \mbox{ for almost every } t \in [0,T].
    \end{eqnarray*}
From the proof of Theorem \ref{maintheorem1} we have  $\tau_{n,t} \rightarrow 0$ for almost every $t \in [0,T]$, a.s.  By the a.s. convergence of $\tilde{X}^{\varepsilon_n}$ to $\tilde{X}$ in $C([0,T]; \R^d)$, we also have $\delta_n \rightarrow 0$. Thus $\gamma_{n,t} = \tau_{n,t} + K \delta_n \rightarrow 0$ for almost every $t \in [0,T]$, a.s.

Take $\bar{\gamma}_{n,t} = \sup_{m \geq n} \gamma_{m,t}$. We then have \begin{eqnarray}
  \label{eqn:derivative-set}
  \frac{d}{dt}\tilde{\xi}^{\varepsilon_n}_t \in h_E(\tilde{X}_t) + \bar{\gamma}_{n,t} \bar{\mathbb{B}}_1 \mbox{ for almost every } t \in [0,T], \forall n \geq 1,
\end{eqnarray}
and $\bar{\gamma}_{n,t} \rightarrow 0$ for almost every $t \in [0,T]$, a.s. Since $\frac{d}{dt}\tilde{\xi}^{\varepsilon_n} \rightarrow \frac{d}{dt}\tilde{X}$ weakly in $L^2{([0,T] \times \tilde{\Omega}; \R^d)}$, a.s., and on account of (\ref{eqn:derivative-set}), by Mazur's lemma (\cite[Lemma ~10.19]{renardy2006introduction}), there exists $\{ Z_n \}_{n \geq 1}$ such that
\[
  Z_n \rightarrow \frac{d}{dt}\tilde{X} \mbox{ in } L^2{([0,T] \times \tilde{\Omega}; \R^d)} \mbox{ as } n \rightarrow \infty,
\]
and
\[
  Z_{n,t} \in co \left( \bigcup_{m \geq n} \left\{ h_E(\tilde{X}_t) + \bar{\gamma}_{m,t} \bar{\mathbb{B}}_1 \right\} \right) \mbox{ for almost every } t \in [0,T].
\]
By passing to a further subsequence, we have $Z_{n,t} \rightarrow \frac{d}{dt}\tilde{X}_t$ for almost every $t \in [0,T]$, a.s. Thus almost surely and for almost every $t \in [0,T]$, we have:
\begin{eqnarray*}
  \frac{d}{dt}\tilde{X}_t & \in & \bigcap_{n \geq 1} co \left( \bigcup_{m \geq n} \left\{ h_E(\tilde{X}_t) + \bar{\gamma}_{m,t} \bar{\mathbb{B}}_1 \right\} \right) \\
  & = & \bigcap_{n \geq 1} co \left\{ h_E(\tilde{X}_t) + \bar{\gamma}_{n,t} \bar{\mathbb{B}}_1 \right\}\\
   & & \quad \mbox{(because $\bar{\gamma}_{n,t}$ neighborhood contains all others for $m \geq n$)}\\
  & = & \bigcap_{n \geq 1} \left( \left\{ h_E(\tilde{X}_t) + \bar{\gamma}_{n,t} \bar{\mathbb{B}}_1 \right\} \right) \\
  & & \quad \mbox{(because $h_E(\tilde{X}_t)$ is already convex and so is its $\bar{\gamma}_{n,t}$ neighborhood)} \\
  & = & h_E(\tilde{X}_t) \quad \mbox{(because $h_E(\tilde{X}_t)$ is also closed)}.
\end{eqnarray*}
By suitably modifying $\frac{d}{dt}\tilde{X}$ on a Lebesgue null set, we establish that $\frac{d}{dt}\tilde{X}_t \in h_E(\tilde{X}_t)$ for all $t \in [0,T]$. Finally, since $X$ and $\tilde{X}$ have the same law, we conclude that, almost surely, $\frac{d}{dt}X_t \in h_E(X_t)$ for all $t \in [0,T]$.

We now argue that any limit point in law is almost surely a Filippov solution to (\ref{odemain}). Let $\delta_n \rightarrow 0$. Along a subsequence,  $X^{\delta_n}$ converges weakly to a limit point $X$ as $\delta_n \rightarrow 0$. There is a further subsequence along which Proposition \ref{p:tightness} and Proposition \ref{characterisation} hold. Imitating the steps of the proof of the first part above along this subsequence, we see that the limit point $X$ is almost surely a Filippov solution to (\ref{odemain}).
\qed

\section{Proof of Proposition \ref{p:spgp} }  \label{pop:spgp}
{A spectral gap estimate is shown in \cite[Theorem 3.1]{HS90}. To convert the estimate in our setting and  to the required $L_\infty$ norm as stated in Proposition \ref{p:spgp}(b) will require ultracontractivity bounds. For this we will need one additional notation. For $1\leq p,q \leq \infty$, write $\| \cdot \|_{(p,q)}$ for the $L_p \rightarrow L_q$ operator norm, with $L_p$ being the space of functions whose $p$-th power is integrable. Our first lemma establishes ultracontractivity.}
\begin{lemma}\label{l:smgp} Let $x \in \R^d, \varepsilon >0$ and $\eta >1$ be as in \eqref{U23}. For $0< t_0 < 1$,  there exists $c_1  > 0$ such that
  \ben \| T_{t_0}^{0,\varepsilon,x}\|_{(1,\infty)} < \exp \left ({\frac{c_1t_0^{-\frac{\eta}{\eta-1}}}{s(\varepsilon)^2}} \right).\label{eqn:hypercontractivity}\een
\end{lemma}

{\em Proof:} Fix $0< t_0 < 1$. The result follows directly from (\ref{U23}) of Assumption (U3) with $a = s(\varepsilon)^2$, \cite[Proposition 7.3.1]{BGL14}, and \cite[Corollary 7.1.4]{BGL14} with  $W(\cdot) = 2U(x, \cdot)/s(\varepsilon)^2$ and the growth function $$\Phi(r) = \frac{C}{s(\varepsilon)^2}(1 + r^{\frac{\eta}{2\eta-1}}), \mbox{ with }  r \in (0,\infty) \mbox{ and } C \equiv C(M, m, \eta).$$

In particular, see the discussion in \cite[p.~363]{BGL14} explaining  the choice of the above growth function $\Phi$ in \cite[eqn.~(7.3.1)]{BGL14}, which satisfies an entropy-energy inequality (\cite[Defn.~7.1.1]{BGL14}) by virtue of (\ref{U23}) with $a = s(\varepsilon)^2$ and \cite[Proposition 7.3.1]{BGL14}. Then \cite[Corollary 7.1.4]{BGL14} yields \eqref{eqn:hypercontractivity}.
\qed

{\em Proof of Proposition \ref{p:spgp}: } From (\ref{U21}),(\ref{U22}),(\ref{U23}) with $a=1$, we may conclude
\bena
 &&  C(x,\varepsilon) := \int_{\R^m} \exp \left\{ - 2\frac{U(x,y)}{s(\varepsilon)^2} \right\} dy < \infty, \label{U24}\\
   &&\norm{\nabla_y U(x,y)}  \rightarrow  \infty \mbox{ as } \norm{y} \rightarrow \infty, \label{U25}\\
 &&  \norm{\nabla_y U(x,\cdot)}^2 - \Delta_y U(x,\cdot)  \mbox{ is bounded below}, \label{U26}\\
   && U(x,y)  \rightarrow  \infty \ \mbox{ as } \norm{y} \rightarrow \infty, \mbox{ uniformly in } x. \label{U27}
  \eena

  So, part (a) follows from the results of Appendix \ref{existenceanduniqueness} along with
  the fact that $\LA^{\varepsilon,x}$ in (\ref{generator}) is a self-adjoint operator on $L^2(\nu^{\varepsilon,x})$ and $\nu^{\varepsilon,x}(\R^m) = 1$.

(b) Using a standard result on  spectral gap (see discussion on \cite[p.~273]{HS90}), we have for all $s \geq t$
\begin{equation}\label{eq6identification}
  \| T_{s-t}^{0,\varepsilon,x}   - \nu^{\varepsilon,x}\|_{(2,2)} \leq e^{-(s-t)\lambda^\varepsilon_2(x)},
\end{equation}
 where $\lambda^\varepsilon_2(x)$ is the second largest eigenvalue of $\LA^{\varepsilon,x}$.
 {Let  $t_0 =\frac{1}{2}, s \geq t + \frac{1}{2}$.}  Using Lemma \ref{l:smgp}
 \begin{eqnarray*}
  \| T_s^{t,\varepsilon,x} -\nu^{\varepsilon,x} \|_{(\infty,\infty)} & =& \| T_{s-t}^{0,\varepsilon,x} -\nu^{\varepsilon,x} \|_{(\infty,\infty)}\\  &\leq & \| T_{s-t}^{0,\varepsilon,x} -\nu^{\varepsilon,x} \|_{(2,\infty)}\\
  &\leq & \| T_{t_0}^{0,\varepsilon,x}\|_{(2,\infty)} \| T_{s-(t+t_0)}^{0,\varepsilon,x} -\nu^{\varepsilon,x} \|_{(2,2)}\\
    &\leq & \| T_{t_0}^{0,\varepsilon,x}\|_{(1,\infty)} \| T_{s-(t+t_0)}^{0,\varepsilon,x} -\nu^{\varepsilon,x} \|_{(2,2)}\\
  &\leq &{e^{\frac{c_12^{\frac{\eta}{\eta-1}}}{s(\varepsilon)^2}}e^{-(s-(t+\frac{1}{2}))\lambda^\varepsilon_2(x)}.}
  \end{eqnarray*}
So for all $f \in C^2_b(\R^m)$, there exists $c_2 >0$ such that
\begin{equation}\label{eq7}
\| T_s^{t,\varepsilon,x} f  - \nu^{\varepsilon,x}(f) \|_{\infty}
\leq  \| f \|_{\infty} e^{\frac{c_2}{s(\varepsilon)^2} -(s-t-\frac{1}{2})\lambda^\varepsilon_2(x)}.
\end{equation}
As (\ref{U24}), (\ref{U25}), (\ref{U26}), (\ref{U27}) hold, from \cite[Theorem 3.1]{HS90}, we obtain that for any $\delta_1 >0$,
$$\lambda^\varepsilon_2(x) \geq \exp\left\{-\frac{V^{(1)}(x) - V^{(2)}(x)+\delta_1}{s(\varepsilon)^2}\right\}$$ for all sufficiently small $\varepsilon$. Using (\ref{u3}), from assumption (U3), and (\ref{eq7}) we have for some $c_3 >0$
\begin{equation}\label{eq8}
\| T_s^{t,\varepsilon,x} f  - \nu^{\varepsilon,x}(f) \|_{\infty}
\leq { \| f \|_{\infty}e^{\frac{c_2}{s(\varepsilon)^2} -(s-t -\frac{1}{2}) \exp(-\frac{\Lambda + \delta_1}{s(\varepsilon)^2})} \leq   \| f \|_{\infty} e^{\frac{c_3}{s(\varepsilon)^2} -(s-t) \exp(-\frac{\Lambda + \delta_1}{s(\varepsilon)^2})}}.
\end{equation}
\qed

\section{Proof of Proposition \ref{p:tightness}} \label{pop:tightness}

 It is easy to obtain fourth moment bounds for
$X^\varepsilon$ from the assumption (B1), this readily implies
tightness, and consequently part (a). Part (b) is a standard
application of Skorohod's Theorem. As indicated earlier the key
nuance in  the Proposition is the topology on $\widetilde{\PA}$. One
of the facts we shall crucially use is that $\widetilde{\PA}$ is
compact and metrizable in this topology. This and other applications to control
theoretic setting are discussed in detail in \cite{ABG12}.

{\em Proof of Proposition \ref{p:tightness}:} (a) Let $0 < \varepsilon <1$ and $0 \leq s \leq t <T$. As $X^{\varepsilon}_t$ solves (\ref{ex1}) we have
 \benas
 \| X^{\varepsilon}_t - X^{\varepsilon}_s \|  &=& \left\| \int_s^t \int_{\R^m}b(X^{\varepsilon}_r,y )\pi^{\varepsilon}_r(dy) dr + \varepsilon^{\alpha}(\eta^\varepsilon_t - \eta^\varepsilon_s) \right\|\\
 &\leq &\| b\|_\infty (t-s) + \varepsilon^{\alpha}\|\eta^\varepsilon_t - \eta^\varepsilon_s\|.
 \eenas
We can then conclude that
$$\E\left[\|X^{\varepsilon}_t - X^{\varepsilon}_s\|^4\right] \leq c_1|t - s|^2$$
 for $0 \leq s < t < T$. By \cite[(12.51) and Theorem 12.3]{B68} we have that the laws of $\{X^{\varepsilon}: \varepsilon \in (0, 1]\}$ are tight  in $\PA(C([0, T]; \mathbb{R}^d))$.
Further, we note that $\widetilde{\PA}$ is compact and metrizable \cite[Section 2.3, Theorem  2.3.1]{ABG12}. This implies the tightness of the laws of $(X^\varepsilon, \pi^\varepsilon)$  in $\PA ( C( [0,T]; \mathbb{R}^d) \times \widetilde{\PA}).$ Hence there exists a sequence $\varepsilon_n \downarrow 0$  such that
$(X^{\varepsilon_n}, \pi^{\varepsilon_n})$ converges weakly to $(X, \pi)$ as $n \rightarrow \infty$.

(b)  Let $\{\varepsilon_n\}_{n \geq 1}$ be the sequence mentioned  in part (a). Using Skorohod's theorem [\cite{B95}, Theorem 2.2.2, p.\ 23], there exists a probability space $(\tilde{\Omega}, \tilde{\mathcal{ F}}, \tilde{P})$ and processes
$(\tilde{X}^{\varepsilon_n}, \tilde{\pi}^{\varepsilon_n}, \tilde{\eta}^{\varepsilon_n}, \tilde{X}, \tilde{\pi}, \widetilde{\eta}) $ such that
\begin{eqnarray*}
{\rm Law\ of}\ (\tilde{X}^{\varepsilon_n}, \tilde{\pi}^{\varepsilon_n}, \tilde{\eta}^{\varepsilon_n})
& = &  {\rm Law\ of}\ (X^{\varepsilon_n}, \pi^{\varepsilon_n}, \eta^{\varepsilon_n}),\\
{\rm Law\ of}\ (\tilde{X}, \tilde{\pi}, \tilde{\eta})
& = &  {\rm Law\ of}\ (X, \pi, \eta),
\end{eqnarray*}
and $(\tilde{X}^{\varepsilon_n} , \tilde{\pi}^{\varepsilon_n}, \tilde{\eta}^{\varepsilon_n}) \to (\tilde{X},  \tilde{\pi}, \tilde{\eta})$ almost surely. Further, using Fatou's lemma followed by Doob's inequality we have
\benas
\E \left[\sup_{0 \leq s \leq T} \| \tilde{X}^{\varepsilon_n}_s - \tilde{X}_s \|^4 \right] & \leq & \liminf_{\delta \rightarrow 0}  \E \left[ \sup_{0 \leq s \leq T} \| \tilde{X}^{\varepsilon_n}_s - \tilde{X}^\delta_s \|^4 \right]\\
    &\leq& 2^4 \left ( \| b \|^4_\infty s^4 + (\varepsilon^\alpha + \delta^\alpha)^4 \E \left[\sup_{0 \leq s \leq T}\| B_s\|^4 \right]\right)\\
    &\leq& c_1\left ( \| b \|^4_\infty  + \E[ \| B_T\| ^4]\right) < \infty.
\eenas
This implies that the family
\[
  \left\{ \sup_{0 \leq s \leq T} \| \tilde{X}^{\varepsilon_n}_s - \tilde{X}_s \|^2: n \geq 1 \right\}
\]
is uniformly integrable. This implies (\ref{eqn:second-moment-convergence}).
  \qed

\section{Proof of Proposition \ref{characterisation}}  \label{pop:characterisation}

{The proof of this proposition consists of many steps. Part (a)
  uses the topology on $\widetilde{\PA}$ and fundamental theorem of
  calculus to choose an appropriate subsequence. Part (b) and Part (c) require some
  technical preparation which we describe in detail first, before proving Proposition \ref{characterisation}.

For Part (b), we prove a second moment estimate in Proposition
\ref{l:ke}. Using this second moment estimate we will be able to
identify the required rate of decay of $s(\varepsilon ) \rightarrow 0$
as $\varepsilon \rightarrow 0.$ This will ensure that the second
moment goes to zero and consequently a further  subsequence goes to zero
almost surely.  Proof of Proposition \ref{l:ke} will require a
gradient estimate for the semigroup of $Z^{t,\varepsilon,x}$ which
satisfies (\ref{associatedsde}). We present that first.}

\begin{lemma} \label{l:ges}
 Recall $\Gamma$ from (\ref{eqn:PW}). There exists $ \varepsilon_0 > 0$, such that for all  $f \in C^2_b(\R^m)$, $s \geq t$, $0 < \varepsilon < \varepsilon_0$,
\ben \label{eqn:gradest} \| \nabla T_{s}^{t,\varepsilon,x}(f)\|_\infty <  \|
f\|_\infty \frac{ e^{\frac{{\Gamma}}{s(\varepsilon)^2}}}{\sqrt{s-t}} .\een
\end{lemma}

{\em Proof:}  Using \eqref{eqn:PW} and \cite[Theorem 3.4]{PW06}  we have  that for any $f \geq 0$ and $f \in C^2_b (\R^m)$,
\[  \| \nabla T_{s}^{t,\varepsilon,x}(f)\|_\infty < \frac{1 + 2 s(\varepsilon)^2}{2 s(\varepsilon)^2\sqrt{s-t}} \exp \left[ \frac{\Gamma }{2s(\varepsilon)^2} \right] \| f\|_\infty. \]
We may choose $\varepsilon_0 >0$ so that $1+1/(2 s(\varepsilon_0)^2) \leq e^{\Gamma/(2 s(\varepsilon_0)^2)}$, and so (\ref{eqn:gradest}) holds.
For any $f \in C^2_b(\R^m)$ the result follows by considering positive and negative parts of $f$.
\qed

{We now present the key second moment estimate.}

\begin{proposition} \label{l:ke} Let $ 0 \leq \Lambda < \infty$ be as in (U3) and let $\delta >0$  be fixed. {There exist $c_1, {c_2} >0$  such that for all $ f \in C^2_b (\mathbb{R}^m)$, for all sufficiently small $\varepsilon >0$,}  {$t \geq 0$, $s > t +1 $,} and  $\kappa >0$,
  \begin{eqnarray} \label{ksmest} \lefteqn{\E \left [ \nu^{\varepsilon,X^\varepsilon_t}(f) -  \frac{1}{\kappa} \int^{t + \kappa}_t {\pi}^\varepsilon_r (f) dr \right ]^2 \,  } \nonumber \\
&& \leq  c_1 \|f\|^2_\infty
\left [{e^{\frac{c_2}{s(\varepsilon)^2} - 2(s-t) \exp\left (-\frac{\Lambda + \delta}{s(\varepsilon)^2} \right)}} + (s-t)\kappa (\kappa + \varepsilon^{2\alpha})e^{\frac{{2\Gamma}}{s(\varepsilon)^2}} + (s-t)^2\left ( {\frac{\varepsilon^2}{\kappa^2}} + \frac{\varepsilon^{-2\alpha +2}}{\kappa} \right) \right]. \nonumber \\
    \end{eqnarray}
\end{proposition}

\begin{remark} \label{spgpsm}
  The first term inside the bracket in (\ref{ksmest}) arises from the
  spectral gap estimate obtained earlier and it specifies the rate at
  which the fast process approaches its stationary measure. The second
  term inside the bracket in (\ref{ksmest}) is from the gradient
  estimate obtained in Lemma \ref{l:ges}. So for both these  terms to go
  to zero, we need to impose a rate of decay to $0$ on
  $s(\varepsilon)$ and use Assumption (S1). The third term
  contains the scaling factor provided by the nonlinear filtering
  equation and here we require $0< \alpha < 1$ for this term to go to $0$.
\end{remark}
{\em Proof of Proposition \ref{l:ke}:} Let $ f \in C^2_b(\R^d), {t \geq 0,  s >  t +1,}  \varepsilon >0, x \in {\mathbb R}^d$  be given.  For $0 <\kappa <1 $, define for notational convenience
\begin{equation}\label{barnu}
\Bar{\nu}^{\varepsilon, \kappa}_t (f)  \ = \ \frac{1}{\kappa} \int^{t + \kappa}_t {\pi}^\varepsilon_r (f) dr.
\end{equation}

By the fundamental theorem of calculus,
\[T^{\varepsilon,t,x}_s f - f = \int^s_t  \LA^{\varepsilon, x}(T^{\varepsilon,t,x}_{u} f) du.\]
We then readily note that
\begin{eqnarray*}
 \nu^{\varepsilon,x}(f) -   \bar{\nu}^{\varepsilon, \kappa}_t (f) =
  \nu^{\varepsilon,x}\left (T^{\varepsilon,t,x}_s f - \int^s_t  \LA^{\varepsilon,x}(T^{\varepsilon,t,x}_{u} f)du \right) -   \bar{\nu}^{\varepsilon, \kappa}_t \left (T^{\varepsilon,t,x}_s f - \int^s_t  \LA^{\varepsilon,x}(T^{\varepsilon,t,x}_u f) du\right).
\end{eqnarray*}
As $f \in C^2_b(\R^d)$ and  $\nu^{\varepsilon,x}$ is  an invariant measure, we have
\ben\label{s2ke}
\nu^{\varepsilon,x}\left ( \int^s_t  \LA^{\varepsilon,x}( T^{\varepsilon,t,x}_{u} f) du\right)= \int^s_t  \nu^{\varepsilon,x}\left ( \LA^{\varepsilon,x}(T^{\varepsilon,t,x}_{u} f)\right) du = 0 .
\een
Using  (\ref{s2ke}), we may rewrite
\begin{eqnarray*}
 \nu^{\varepsilon,x}(f) -   \bar{\nu}^{\varepsilon, \kappa}_t (f) &=&
  \nu^{\varepsilon,x}\left (T^{\varepsilon,t,x}_s f  \right) - \bar{\nu}^{\varepsilon, \kappa}_t \left (T^{\varepsilon,t,x}_s f \right ) + \int^s_t  \bar{\nu}^{\varepsilon, \kappa}_t \left ( \LA^{\varepsilon,x}(T^{\varepsilon,t,x}_u f) \right) du.
\end{eqnarray*}
As both $\nu^{\varepsilon,x}$ and $\bar{\nu}^{\varepsilon, \kappa}_t$ are
probability measures, we may add and subtract the constant term
$\nu^{\varepsilon,x}(f)$ in the first two terms above. Using the definition of $\bar{\nu}^{\varepsilon, \kappa}$ in third term above, we have
\begin{eqnarray}
\lefteqn{\nu^{\varepsilon,x}(f) -   \bar{\nu}^{\varepsilon, \kappa}_t (f) }\nonumber \\
&=&  \nu^{\varepsilon,x}\left (T^{\varepsilon,t,x}_s f  -\nu^{\varepsilon,x}(f) \right) -   \bar{\nu}^{\varepsilon, \kappa}_t \left (T^{\varepsilon,t,x}_s f   -\nu^{\varepsilon,x}(f) \right )  + \int^s_t  \bar{\nu}^{\varepsilon, \kappa}_t \left ( {\LA^{\varepsilon, x}}( T^{\varepsilon,t,x}_u f) \right) du\nonumber \\
  &=&  \nu^{\varepsilon,x}\left (T^{\varepsilon,t,x}_s f  -\nu^{\varepsilon,x}(f) \right) -   \bar{\nu}^{\varepsilon, \kappa}_t \left (T^{\varepsilon,t,x}_s f   -\nu^{\varepsilon,x}(f) \right ) + \int^s_t \left( \frac{1}{\kappa} \int^{t + \kappa}_t \pi^\varepsilon_r ( \LA^{\varepsilon,x}(T^{\varepsilon,t,x}_u f) ) dr\,\right)  du\nonumber \\
  &=&  \nu^{\varepsilon,x}\left (T^{\varepsilon,t,x}_s f  -\nu^{\varepsilon,x}(f) \right) -   \bar{\nu}^{\varepsilon, \kappa}_t \left (T^{\varepsilon,t,x}_s f   -\nu^{\varepsilon,x}(f) \right )\nonumber \\
&& +\int^s_t \left ( \frac{1}{\kappa} \int^{t + \kappa}_t \pi^\varepsilon_r ( \LA^{\varepsilon,x} (T^{\varepsilon,t,x}_u f) - \LA^{\varepsilon, X^\varepsilon_r} (T^{\varepsilon,t,x}_u f) ) dr \right) du \nonumber \\
&&+ \int^s_t  \left ( \frac{1}{\kappa} \int^{t + \kappa}_t \pi^\varepsilon_r \left ( \LA^{\varepsilon,X^\varepsilon_r}(T^{\varepsilon,t,x}_u f) \right) dr  \right) du. \label{in-filt}
\end{eqnarray}
Now, the measure valued process $\pi^\varepsilon$ is the unique solution to the (Fujisaki-Kallianpur-Kunita) nonlinear
filtering equation
\begin{equation}\label{KushnerStratonovicha}
 \pi^\varepsilon_t(f) \ = f(y_0) + \ \frac{1}{\varepsilon}\int_0^t \pi^\varepsilon_r(
{\LA^{\varepsilon , X^\varepsilon_r}(f))} dr
+ \varepsilon^{{- \alpha}} \int_0^t\langle \pi^\varepsilon_r(f b(X^\varepsilon_r, \cdot)) -
\pi^\varepsilon_r (f) \pi^\varepsilon_r(b(X^\varepsilon_r, \cdot)), d \tilde{B}_r \rangle,
\end{equation}
where {$\tilde{B}$} is a standard Brownian motion; see Proposition \ref{fkksde} in Appendix\footnote{Proposition \ref{fkksde}, presented in Appendix, is a more general nonlinear filtering equation and could be of independent interest.}.  Using the definition of $\LA^{\varepsilon, x}$ from (\ref{generator}) and the FKK equation  (\ref{KushnerStratonovicha}) in (\ref{in-filt}) we have,
\begin{eqnarray*}
\lefteqn{\nu^{\varepsilon,x}(f) -   \bar{\nu}^{\varepsilon, \kappa}_t (f) }\\
&=&  \nu^{\varepsilon,x}\left (T^{\varepsilon,t,x}_s f  -\nu^{\varepsilon,x}(f) \right) -   \bar{\nu}^{\varepsilon, \kappa}_t \left (T^{\varepsilon,t,x}_s f   -\nu^{\varepsilon,x}(f) \right )\\
&&+ \int^s_t \left ( \frac{1}{\kappa} \int^{t + \kappa}_t \pi^\varepsilon_r ( \langle \nabla_yU(X^{\varepsilon}_r, \cdot) - \nabla_yU(x, \cdot),
\nabla_y T^{\varepsilon,t,x}_u f\rangle ) dr \right) du \nonumber\\&&   +{\frac{\varepsilon}{\kappa}} \int_s^ t\left ( \pi_{t+\kappa}( T^{\varepsilon,t,x}_u f) -\pi^\varepsilon_t(T^{\varepsilon,t,x}_u f) \right) du \nonumber\\&&
- \int_t^s \left (\frac{\varepsilon^{-{ \alpha} +1}}{\kappa} \int_t^{t+\kappa}\langle \pi^\varepsilon_r(T^{\varepsilon,t,x}_uf b(X^\varepsilon_r, \cdot)) - \pi^\varepsilon_r (T^{\varepsilon,t,x}_uf) \pi^\varepsilon_r(b(X^\varepsilon_r, \cdot)), d \tilde{B}_r \rangle \right) du. \nonumber
\end{eqnarray*}
We shall now replace $x$ in above by $X^{\varepsilon}_t$. To do this one needs to be careful only in the last term. Here we observe that as $t>0$ is fixed, $t \leq r \leq t+\kappa,$ using definition of $\pi^\varepsilon$ and the stochastic integral, we may replace $x$ by $X^\varepsilon_t$. For the other terms, the substitution is trivial. So we have,
\begin{eqnarray}  \label{s1ke}
\lefteqn{\nu^{\varepsilon,{X^{\varepsilon}_t}}(f) -   \bar{\nu}^{\varepsilon, \kappa}_t (f) }\nonumber \\
&=&  \nu^{\varepsilon,{X^{\varepsilon}_t}}\left (T^{\varepsilon,t,{X^{\varepsilon}_t}}_s f  -\nu^{\varepsilon,{X^{\varepsilon}_t}}(f) \right) -   \bar{\nu}^{\varepsilon, \kappa}_t \left (T^{\varepsilon,t,{X^{\varepsilon}_t}}_s f   -\nu^{\varepsilon,{X^{\varepsilon}_t}}(f) \right )\nonumber \\
&&+ \int^s_t \left ( \frac{1}{\kappa} \int^{t + \kappa}_t \pi^\varepsilon_r ( \langle \nabla_yU(X^{\varepsilon}_r, \cdot) - \nabla_yU(X^{\varepsilon}_t, \cdot),
\nabla_y T^{\varepsilon,t,{X^{\varepsilon}_t}}_u f\rangle ) dr \right) du \nonumber\nonumber \\&&   +  { \frac{\varepsilon}{\kappa}} \int_s^ t\left ( \pi_{t+\kappa}( T^{\varepsilon,t,{X^{\varepsilon}_t}}_u f) -\pi^\varepsilon_t(T^{\varepsilon,t,{X^{\varepsilon}_t}}_uf) \right)du \nonumber\nonumber \\&&
 -\int_t^s \left (\frac{\varepsilon^{-{ \alpha} +1}}{\kappa} \int_t^{t+\kappa}\langle \pi^\varepsilon_r(T^{\varepsilon,t,{X^{\varepsilon}_t}}_uf b(X^\varepsilon_r, \cdot)) - \pi^\varepsilon_r (T^{\varepsilon,t,{X^{\varepsilon}_t}}_uf) \pi^\varepsilon_r(b(X^\varepsilon_r, \cdot)), d \tilde{B}_r \rangle \right) du\nonumber \\
&=:& I + II + III - IV.
\end{eqnarray}
So, \begin{equation}\label{po2b}
\E \left [\nu^{\varepsilon,{X^{\varepsilon}_t}}(f) -   \bar{\nu}^{\varepsilon, \kappa}_t (f) \right]^2 \leq 16 ~ \E \left[ I^2 + II^2 + III^2 +IV^2 \right].
\end{equation}
For the first term in (\ref{s1ke}), i.e., I,  by  Proposition \ref{p:spgp} we have that for sufficiently small $\varepsilon >0$
\bena \label{bI}
\lefteqn{ \E [I^2] \leq  \E\left[\mid \nu^{\varepsilon,X^{\varepsilon}_t}\left (T^{\varepsilon,t,X^{\varepsilon}_t}_s f  -\nu^{\varepsilon,X^{\varepsilon}_t}(f) \right) -   \bar{\nu}^{\varepsilon, \kappa}_t \left (T^{\varepsilon,t,X^{\varepsilon}_t}_s f   -\nu^{\varepsilon,X^{\varepsilon}_t}(f) \right ) \mid^2\right]}\nonumber\\
&\leq& 4 \left ( \E\left[\mid \nu^{\varepsilon,X^{\varepsilon}_t}\left (T^{\varepsilon,t,X^{\varepsilon}_t}_s f  -\nu^{\varepsilon,X^{\varepsilon}_t}(f) \right)\mid^2\right] +    \E\left[\mid  \bar{\nu}^{\varepsilon, \kappa}_t \left (T^{\varepsilon,t,X^{\varepsilon}_t}_s f   -\nu^{\varepsilon,X^{\varepsilon}_t}(f) \right )  \mid^2\right] \right) \nonumber\\
&\leq & 8 \sup_{x \in \R^d} \parallel T^{\varepsilon,t,x}_s f  -\nu^{\varepsilon,x}(f) \parallel^2_\infty\nonumber\\
&\leq&  c_3 {e^{\frac{c_4}{s(\varepsilon)^2}- 2(s-t) \exp\left (-\frac{\Lambda + \delta}{s(\varepsilon)^2} \right)}} \|f\|^2_\infty.
\eena
For the second term in (\ref{s1ke}), i.e., II, using (U1) and (\ref{eqn:gradest}), we have that
\begin{eqnarray}\label{bII}
  \lefteqn{\E[II^2] \leq \E\left| \int^s_t \left (  \frac{1}{\kappa} \int^{t + \kappa}_t \pi^\varepsilon_r ( \langle \nabla_yU(X^{\varepsilon}_r, \cdot) - \nabla_yU(X^\varepsilon_t, \cdot) , \nabla T^{\varepsilon,t,X^\varepsilon_t}_u f\rangle ) dr \right) du \right|^2} \nonumber \\ &\leq&
{c_4 K_2}  (s-t)^2 \sup_{u \in (s,t)} \| \nabla T^{\varepsilon,t,X^\varepsilon_t}_u(f)  \|^2_\infty  \E \left[ \frac{1}{\kappa} \int^{t + \kappa}_t \| X^\varepsilon_r -X^\varepsilon_t \|^2 dr \right]\nonumber\\
&\leq& c_5   e^{\frac{{2\Gamma}}{s(\varepsilon)^2}}(s-t){\|f\|^2_\infty} \frac{1}{\kappa} \int^{t + \kappa}_t  \E \left [ \kappa^2 \| b \|^2_\infty + \varepsilon^{2\alpha}  \| B_{r} - B_t \|^2 \right]dr \nonumber\\
&\leq& c_5   e^{\frac{{2\Gamma}}{s(\varepsilon)^2}}(s-t){\|f\|^2_\infty}  \frac{1}{\kappa} \int^{t + \kappa}_t  \left [ \kappa^2 \| b \|^2_\infty + \varepsilon^{2\alpha}  \E[\|B_{r} - B_t\|]^2 \right]  dr\nonumber\\
&\leq& c_6 e^{\frac{{2\Gamma}}{s(\varepsilon)^2}} (s-t) \kappa \left (\kappa + \varepsilon^{2\alpha}  \right) {\|f\|^2_\infty }.
\end{eqnarray}
For the third term in (\ref{s1ke}), i.e., III, as $\pi^\varepsilon$ is a probability measure,
we have by triangle inequality and the semigroup property,
 {
\bena \label{bIII}
 \E [III^2] \leq \E\left[\frac{\varepsilon^2}{\kappa^2} \left| \int_t^s \left( \pi_{t+\kappa}( T^{\varepsilon,t,x}_u f) -\pi^\varepsilon_t(T^{\varepsilon,t,x}_uf) \right) du \right|^2 \right] &\leq & c_7\frac{\varepsilon^2}{\kappa^2}\,  (s-t)^2 \|f \|^2_\infty .
\eena}
For the fourth term in (\ref{s1ke}), i.e., IV, using Jensen's inequality and a standard second moment estimate, we have
\begin{eqnarray}\label{bIV}
\E [IV^2] &\leq & E\left(\int_t^s \left (\frac{\varepsilon^{-{ \alpha} +1}}{\kappa} \int_t^{t+\kappa}\langle \pi^\varepsilon_r(T^{\varepsilon,t,x}_uf b(X^\varepsilon_r, \cdot)) -\pi^\varepsilon_r (T^{\varepsilon,t,x}_uf) \pi^\varepsilon_r(b(X^\varepsilon_r, \cdot)), d \tilde{B}_r \rangle \right) du \right)^2\nonumber\\
  &\leq& {c_8} (s-t)^2 \frac{\varepsilon^{-2\alpha +2}}{\kappa}  \parallel f \parallel^2_\infty \parallel b \parallel^2_\infty = {c_{9}}  (s-t)^2 \frac{\varepsilon^{-2\alpha +2}}{\kappa}
{\|f\|^2_\infty} .
\end{eqnarray}
So from \eqref{barnu}, \eqref{po2b}, \eqref{bI},\eqref{bII},\eqref{bIII},\eqref{bIV} we have  the result.
\qed

We are now ready to prove the main result of this section.

{\em Proof of Proposition \ref{characterisation}:} Let $f \in {\mathcal C^2_0(\R^m)}.$ Let $\varepsilon_n >0, \tilde{X}, \tilde{\pi}$ be as constructed in Proposition \ref{p:tightness}. Recall that $\gamma = \min\{1-\alpha, \frac{1}{2}\}$.

(a) By the topology of $\tilde{\PA}$, we have for any $\eta >0$
\[
  \lim_{n \to \infty} \int^{t + \eta}_t \tilde{\pi}^{\varepsilon_n}_r
(f) dr \ = \ \int_t^{t + \eta} \tilde{\pi}_r (f)dr, \ \mbox{a.s.}
\]
By \cite[Proposition~7.5.7]{AS08}, we have
\[
  \lim_{\eta \to 0}   \lim_{n \to \infty} \frac{1}{\eta} \int^{t + \eta}_t \tilde{\pi}^{\varepsilon_n}_r
(f) dr \ = \ \tilde{\pi}_t (f), \mbox{ almost every } t \in [0,T], \ \mbox{a.s.}
\]
Hence for each $k \in \N$, there exists a $\eta_k > 0$ and $\varepsilon(k, \eta_k) > 0$ such that
\ben \label{eqn:0overk}
\Big| \frac{1}{\eta_k} \int^{t + \eta_k}_t \tilde{\pi}^{\varepsilon_n}_r (f) dr - \tilde{\pi}_t (f)
\Big| < \frac{1}{k}, \een $\ \forall \varepsilon_n \leq \varepsilon(k, \eta_k),  \mbox{ almost every } t \in [0,T], \ \mbox{a.s.},$ and furthermore, $\eta_k \to 0$ as $k \to \infty$.

For each $k \geq 1$, choose $n_k$ sufficiently large so that both $\varepsilon_{n_k} \leq \varepsilon(k,\eta_k)$ and $\varepsilon_{n_k}^{\gamma} \leq \eta_k$. Then by construction we have the following:
\ben \label{eqn:1overk}
\Big| \frac{1}{\varepsilon_{n_k}^{\gamma}} \int^{t + \varepsilon_{n_k}^{\gamma}}_t \tilde{\pi}^{\varepsilon_{n_k}}_r (f) dr - \tilde{\pi}_t (f)
\Big| < \frac{1}{k}, \een $\ \forall k  \geq 1,  \mbox{ almost every } t \in [0,T], \ \mbox{a.s.}$
The result follows.

(b) Proposition \ref{p:tightness}(b) implies that
\ben\label{tnti}
\tilde{\E}  \left| \nu^{\varepsilon_n,\tilde{X}^{\varepsilon_n}_t}(f) -   \frac{1}{\varepsilon_n^{\gamma}} \int^{t + \varepsilon_n^{\gamma}}_t\tilde{\pi}^{\varepsilon_{n}}_r (f) dr \right|^2 = \E \left| \nu^{\varepsilon_n,X^{\varepsilon_n}_t}(f) -    \frac{1}{\varepsilon_n^{\gamma}} \int^{t + \varepsilon_n^{\gamma}}_t{\pi}^{\varepsilon_{n}}_r (f) dr  \right|^2.\een
Using Proposition \ref{l:ke}, with $\kappa: = \varepsilon_n^{\gamma} $ with $\gamma = \min \{1-\alpha, \frac{1}{2}\}$ and $\delta $ to be chosen soon, we have for sufficiently large $n$
  \begin{eqnarray*}  \lefteqn{\E \left| \nu^{\varepsilon_n,X^{\varepsilon_n}_t}(f) -    \frac{1}{\varepsilon_n^{\gamma}} \int^{t + \varepsilon_n^{\gamma}}_t{\pi}^{\varepsilon_{n}}_r (f) dr  \right|^2}  \\
  &\leq& c_1 {\|f\|^2_\infty} \left [{ e^{ \frac{c_1}{s(\varepsilon)^2}- 2(s-t) \exp\left (-\frac{\Lambda + \delta}{s(\varepsilon_n)^2} \right)}} + (s-t)\varepsilon_n^{\gamma}(\varepsilon_n^{\gamma} + \varepsilon_n^{2\alpha})e^{\frac{2\Gamma}{s(\varepsilon_n^2)}} +(s-t)^2\left ( { \varepsilon_n^{2-2\gamma}} +  \varepsilon_n^{2-2\alpha -\gamma} \right) \right],
   \end{eqnarray*}
  for all $t \geq 0 $ and $s > t + 1 .$ Substituting in the above   $s = t + \varepsilon_n^{-\theta}$, with $\theta >0$ to be chosen soon and  ${s(\varepsilon_n)^2  \geq \frac{C}{\ln (1 +\frac{1}{\varepsilon_n})}}$ we have
\begin{eqnarray*}
   \lefteqn{ \E \left| \nu^{\varepsilon_n,X^{\varepsilon_n}_t}(f) -   \frac{1}{\varepsilon_n^{\gamma}} \int^{t + \varepsilon_n^{\gamma}}_t{\pi}^{\varepsilon_{n}}_r (f) dr  \right|^2} \\
   &\leq& c_1 \|f\|^2_\infty\left [  {e^{c_3 \ln(1 + \frac{1}{\varepsilon_n}) -2 \varepsilon_n^{-\theta + \frac{\Lambda+ \delta}{C}} (1+ {\varepsilon_n})^{-\frac{\Lambda + \delta}{C}}}} +  \varepsilon_n^{-\theta + \gamma}( \varepsilon_n^{\gamma} + \varepsilon_n^{2\alpha} )\varepsilon_n^{-\frac{{2\Gamma}}{C}} (1+\varepsilon_n)^{\frac{{2\Gamma}}{C}} \right.  \\ && \left . \hspace{2in} + \varepsilon_n^{2-2\gamma- 2 \theta} +  \varepsilon_n^{2-2\alpha -\gamma- 2\theta}  \right]\\
   &\leq& {c_2 \|f\|^2_\infty\left [  {e^{c_4 (-\ln(\varepsilon_n) -2 \varepsilon_n^{-\theta + \frac{\Lambda+ \delta}{C}})}} +  \varepsilon_n^{-\theta + 2\gamma-\frac{{2\Gamma}}{C}} +   \varepsilon_n^{-\theta + \gamma + 2\alpha -\frac{{2\Gamma}}{C}} + \varepsilon_n^{2-2\gamma- 2 \theta} +  \varepsilon_n^{2-2\alpha -\gamma- 2\theta}  \right].}\end{eqnarray*}
{If we can choose  $\delta$ and $\theta$ such  that \begin{equation}\label{chothde1}\frac{\Lambda + \delta}{C}< \theta < \min \{   2\gamma -\frac{2\Gamma}{C},\gamma + 2\alpha-\frac{2\Gamma}{C},  1-\gamma, 1-\alpha -\frac{\gamma}{2} \}\end{equation}
then this would imply
  \ben \nonumber
   \lim_{n \rightarrow \infty}    \E \left| \nu^{\varepsilon_n,X^{\varepsilon_n}_t}(f) -  \frac{1}{\varepsilon_n^{\gamma}} \int^{t + \varepsilon_n^{\gamma}}_t{\pi}^{\varepsilon_{n}}_r (f) dr  \right|^2 =0.  \label{il3}
   \een
   which would then imply, along with (\ref{tnti}) that  there exists a further subsequence $\{\varepsilon_{n_k}\}_{k \geq 0}$ such that (\ref{l3}) holds.

   {So to complete the proof we need to find  small enough $\theta >0 , \delta >0$  so that (\ref{chothde1}) holds. This will be possible if
  \begin{align*}
    \frac{\Lambda}{C} < 2\gamma-\frac{{2\Gamma}}{C},\,\,\,\,\,\,\,\,\, \frac{\Lambda}{C} < \gamma + 2\alpha -\frac{{2\Gamma}}{C}, \,\,\,\,\,\,\,\,\, \frac{\Lambda}{C} < 1-\gamma,\mbox{ and }\,\,\,\,\,\,\,\,\,\frac{\Lambda}{C} < 1-\alpha -\frac{\gamma}{2}.    \end{align*}
  The above will be true if a
     \begin{align} \label{chothde}
     C \geq \max\left \{ \frac{\Lambda +2\Gamma}{2\gamma},\,\,\,\,\,\,\,\,\, \frac{\Lambda + 2\Gamma}{\gamma + 2\alpha}, \,\,\,\,\,\,\,\,\, \frac{\Lambda}{1-\gamma},\mbox{ and }\,\,\,\,\,\,\,\,\,\frac{\Lambda}{1-\alpha -\frac{\gamma}{2}} \right\}.    \end{align}
As $\gamma = \min \{1-\alpha, \frac{1}{2}\},$ we have  for  $0 < \alpha < 1$ that $$1-\gamma \geq \frac{1}{2}\,\,\,\,\,\,\,\,\, \mbox{ and } \,\,\,\,\,\,\,\,\,\frac{1-\alpha}{2} \leq \gamma \leq 1-\alpha. \,\,\,\,\,\,\,\,\,$$ So, (\ref{chothde})  will be true if        \begin{align} \label{chothde2}
     C > \max\left \{ \frac{\Lambda +2\Gamma}{1-\alpha},\,\,\,\,\,\,\,\,\, \frac{2\Lambda + 4\Gamma}{1+3\alpha},\,\,\,\,\,\,\,\,\, 2\Lambda,  \mbox{ and }\,\,\,\,\,\,\,\,\,\frac{2\Lambda}{1-\alpha}\right \}.    \end{align}
In (S1) we require $C > \frac{2(\Lambda + 2 \Gamma)}{1-\alpha},$ so (\ref{chothde2}) is true.
}

For part (c) we will need to understand how to characterize  weak limit points of $\nu^{\varepsilon_n, x^{\varepsilon_n}}$ when $\varepsilon_n \rightarrow 0$ and for deterministic $x^{\varepsilon_n} \rightarrow x$. For part (d), under (U4), we will need to verify that the above sequence of measures obeys Laplace's principle. We present these results about deterministic sequence of invariant measures in Lemma \ref{l:ieu} of Appendix \ref{appendix:Laplace}. We now use the result in Lemma \ref{l:ieu} to finish the proof.

(c) Let $\varepsilon_{n_k}$ be a subsequence along which (\ref{l3a}) and (\ref{l3}) hold. Using Proposition \ref{p:tightness}(b), there is a null set $N$ such that $\tilde{X}^{\varepsilon_{n_k}}_t \rightarrow \tilde{X}_t\, \mbox{ as } k \rightarrow  \infty $, (\ref{l3a}) and (\ref{l3}) hold for all $\omega \in N^c$ for almost every $t \in [0,T].$

For a fixed $\omega \in N^c$, using Lemma \ref{l:ieu} of Appendix \ref{appendix:Laplace} with $x^{\varepsilon_{n_k}}=\tilde{X}^{\varepsilon_{n_k}}_t(\omega)$ and $x=\tilde{X}_t(\omega)$, we conclude that $\nu^{\varepsilon_{n_k}, \tilde{X}^{\varepsilon_{n_k}}_t(\omega)}, k \geq 1$ are tight. {Let $\varepsilon_{n_{k_l}}$ be a subsequence along which $\nu^{\varepsilon_{n_{k_l}}, \tilde{X}^{\varepsilon_{n_{k_l}}}_t(\omega)}$ converges weakly to a measure (again by Lemma \ref{l:ieu}) supported on $\arg \min U(\tilde{X}_t(\omega), \cdot)$, which by (U1) is a finite set. Let us denote this measure by $\nu_t^{0,\tilde{X}_t}$. Further, from (\ref{l3a}) and (\ref{l3}), $\nu^{\varepsilon_{n_k}, \tilde{X}^{\varepsilon_{n_k}}_t(\omega)}(f)$ converges to $\tilde{\pi}_t(f)$ for $f \in C_0^2(\R^m)$ and for almost every $t \in [0,T]$. It is now standard to see that $\tilde{\pi}_t(f) = \nu_t^{0,\tilde{X}_t}(f) $ for all $f \in C_0^2(\R^m)$ and for almost every $t \in [0,T]$, and consequently that $\tilde{\pi}_t = \nu_t^{0, \tilde{X}_t}$ for almost every $t \in [0,T]$.}

(d) Follows immediately from (c) and Lemma \ref{l:ieu}(b).
\qed

\begin{appendices}
\section{Existence, Uniqueness, and Gradient Estimates} \label{existenceanduniqueness}

In this section we show that the coupled system (\ref{ex}) and (\ref{wye}) has a unique strong solution. We begin with a technical lemma.

\begin{lemma} \label{l:U24}
Under (U1), (\ref{U21}) and (\ref{U22}) in assumption (U2)  there is $K_4 >0$ and $R^\prime \geq R$ such that
\begin{equation}   \langle \nabla_y U(x, y) , y \rangle > K_4 \|y\|^2, \quad \| y \| > R^\prime. \label{U24aa}\end{equation}
Also, there exists a nonnegative continuous function $g: (0,\infty ) \rightarrow (0,\infty)$ such that
  \begin{eqnarray}
&&  \sup_{z,y \in \R^m: \| z-y \| = r} - \frac{1}{r} \langle \nabla_yU(x,z) - \nabla_yU(x,y), z- y \rangle
  \leq  g(r), \mbox{ for all } r > 0,\nonumber\\
&&\mbox{ with }  \Gamma:= \int_0^\infty g(s) ds < \infty.   \label{eqn:PWa}
  \end{eqnarray}
   \end{lemma}
   {\em Proof :}
We proceed as follows. Let $\mathbb{B}_a$ denote the closed ball of radius $a$ centred at the origin.
For any $y$ with $||y|| > R$, writing $\nabla_y U(x,y) = \nabla_y U(x,0) + \int_0^1 D^2_y U(x,ty) y ~dt$, we have
\begin{eqnarray*}
\langle \nabla_y U(x,y), y \rangle & = & \langle \nabla_y U(x,0), y \rangle + \int_0^1 \langle
y, D_y^2U(x,ty) y \rangle ~dt \\
& \geq & -M||y|| - \int_0^{R/||y||} M||y||^2 dt +  \int_{R/||y||}^1  K_3||y||^2 dt  \\
& = & -M||y|| - MR ||y|| + K_3 ||y||^2 (1 - R/||y||) \\
& = & ( K_3 ||y||^2 - (M + MR + K_3 R)) ||y|| \\
& \geq & K_4 ||y||^2
\end{eqnarray*}
for any $K_4 > 0$ and $\|y\| \geq R^\prime \geq R + M(1+R)/K_3 + K_4 / K_3$. In the second inequality above, we have used \eqref{U22} for the line segment joining $0$ to $y$ that lies within $\mathbb{B}_R$ and \eqref{U21} for the remaining line segment. This establishes \eqref{U24aa}.

Next, for any $y,z \in \R^m$, define $t_0(y,z)$ to be the fractional length of the line segment joining $y$ to $z$ that is within $\mathbb{B}_R$. Take $R_1 = R(1+2Mm/K_3)$. With $r = ||y-z||$, we can write
\begin{eqnarray}
 \lefteqn{ \frac{1}{r} \langle \nabla_yU(x,z) - \nabla_yU(x,y), z- y \rangle } \nonumber \\
 & = & \frac{1}{r} \int_0^1 \langle (z-y), D_y^2U(x,y+t(z-y)) (z-y) \rangle ~dt \nonumber \\
 & = & \frac{1}{r} \int_0^1 \langle (z-y), D_y^2U(x,y+t(z-y)) (z-y) \rangle ~1_{\mathbb{B}_R}(y+t(z-y)) ~ dt \nonumber \\ && \hspace{1in} + \frac{1}{r} \int_0^1 \langle (z-y), D_y^2U(x,y+t(z-y)) (z-y) \rangle ~ 1_{\mathbb{B}_R^c}(y+t(z-y))  ~ dt \nonumber \\
 & \geq & -\frac{1}{r} t_0(y,z) Mmr^2 + \frac{1}{r}(1 - t_0(y,z)) K_3 r^2 \label{eqn:three-cases} \\
 & = & r (K_3 - t_0(y,z) (Mm+K_3)) \nonumber \\
 \label{eqn:three-cases-2}
 & \geq & \begin{cases}
 -Mmr & \mbox{if } y,z \in \mathbb{B}_{R_1} \\
 0 & \mbox{otherwise}.
 \end{cases}
\end{eqnarray}
The inequality in (\ref{eqn:three-cases}) follows because:
\begin{itemize}
 \item[(a)] from (\ref{U22}), on account of $||D^2_yU(x,y+t(z-y))|| \leq M$ when $y+t(z-y) \in \mathbb{B}_R$, we easily obtain the simple inequality $\langle (z-y), D_y^2(x, y+t(z-y))(z-y) \rangle \geq -Mmr^2$ using which the first term is obtained; and
 \item[(b)] from (\ref{U21}), $\langle (z-y), D_y^2(x, y+t(z-y))(z-y) \rangle \geq K_3 r^2$ when $y+t(z-y)$ is outside $\mathbb{B}_R$.
\end{itemize}
The inequality in (\ref{eqn:three-cases-2}) follows from the easily verifiable fact
\begin{equation}
t_0(y,z) \leq 2R/(R+R_1) = K_3/(K_3 + Mm)
\end{equation}
for every $y,z$ such that one of them is outside $\mathbb{B}_{R_1}$.

From (\ref{eqn:three-cases-2}), it is clear that we may take $g(\cdot)$ to be any continuous function that dominates the function $Mmr \cdot {\bf 1}\{ r \leq 2R_1\}$, and there is at least one such $g(\cdot)$ that satisfies $\int_0^{\infty} g(s) ds < \infty$.
\qed

Given Brownian motion $B_t$ on $\R^d$ and an independent Brownian
motion $W_t$ on a filtered probability space $(\Omega, {\cal F},
 \P)$, a strong solution to the coupled system (\ref{ex}) and
 (\ref{wye}) is a continuous process $(X^{\varepsilon}_t,
 Y^{\varepsilon}_t)$ that is adapted to the complete filtration
 generated by $B,W$, and satisfies (\ref{ex}) and (\ref{wye}). We
 say that strong uniqueness holds for the coupled system (\ref{ex})
 and (\ref{wye}) if whenever $(X^{\varepsilon}_t, Y^{\varepsilon}_t)$
 and $(\tilde{X}^{\varepsilon}_t, \tilde{Y}^{\varepsilon}_t)$ are two
 strong solutions of the coupled system (\ref{ex}) and (\ref{wye})
 with the common initial condition $x_0,y_0$, then
 $\P((X^{\varepsilon}_t, Y^{\varepsilon}_t) =
 (\tilde{X}^{\varepsilon}_t, \tilde{Y}^{\varepsilon}_t) \mbox{ for
   all } t \geq 0) = 1.$

\begin{lemma} \label{lemma:eunexp} Assume (B1), (U1) and (U2). Let $\varepsilon >0, 0 < \alpha < 1$ and $s(\varepsilon) >0$ be given. The coupled system given by (\ref{ex}) and (\ref{wye}) has  a unique strong solution.
\end{lemma}

   {\em Proof} By assumptions (B1) and (U1), we know that $b: \R^d
   \times \R^m { \to \mathbb{R}^d}$ and
$\nabla_yU: \R^d \times \R^m { \to \mathbb{R}^m}$ are locally
   Lipschitz functions. By \cite[page 178 Theorem 3.1]{IW89} there exist
   a unique strong solution, $({X}^{\varepsilon}_t,
   {Y}^{\varepsilon}_t)_{0 \leq t < \zeta}$ where $$ \zeta =  \inf \{ t \geq 0: \|{X}^{\varepsilon}_t\|^2 + \|{Y}^{\varepsilon}_t\|^2 = \infty \}.$$
We will now establish nonexplosiveness of the process. Let $f : \R^d \times \R^m \rightarrow [0, \infty)$ be given by $f(x,y) = ~\|x \|^2 + \| y \|^2$. Let $$\sigma_n =\inf \{ t \geq 0: \| {X}^{\varepsilon}_t\|^2 + \|{Y}^{\varepsilon}_t\|^2 = n\}.$$ Clearly $\sigma_n \leq \zeta$ almost surely for all $n \geq 1$. Let $t >0$ be given. Applying Ito's formula at time $\sigma_n \wedge t$,
       we obtain that
       \begin{eqnarray}
         \E[f(X^{\varepsilon}_{\sigma_n \wedge t},Y^{\varepsilon}_{\sigma_n \wedge t})] &=& f(x_0,y_0) +  \E\int_0^{\sigma_n \wedge t}2 \left (\langle  X^{\varepsilon}_r , b( X^{\varepsilon}_r,  Y^{\varepsilon}_r) \rangle - \langle  Y^{\varepsilon}_r , \frac{1}{\varepsilon}\nabla_yU( X^{\varepsilon}_r,  Y^{\varepsilon}_r) \rangle \right) dr \nonumber\\ &&+  (m s(\varepsilon)^2/\varepsilon + d \varepsilon^{2\alpha}) \E(\sigma_n \wedge t). \label{nesp}
       \end{eqnarray}
       Using the fact that $b$ is bounded from assumption (B1) we have, for $r >0$,
      \begin{equation} \langle  X^{\varepsilon}_r , b( X^{\varepsilon}_r,  Y^{\varepsilon}_r) \rangle \leq c_1 d \| X^{\varepsilon}_r\|  \|b \|_\infty. \label{xbound}\end{equation}
 Using (\ref{U22}) from  assumption (U2) and (\ref{U24aa}) derived above we have, for $r >0$,
 \begin{equation} - \langle  Y^{\varepsilon}_r , \frac{1}{\varepsilon}\nabla_yU( X^{\varepsilon}_r,  Y^{\varepsilon}_r) \rangle < \begin{cases} c_2 (M,R) &\mbox{if } { \|Y_r\| \leq R}\cr
     0 &\mbox{if } { \|Y_r\| > R}.
     \end{cases}\label{ybound} \end{equation}
 Substituting (\ref{xbound}) and (\ref{ybound}) in (\ref{nesp}) we have
     \begin{eqnarray}
       \E(f(X^{\varepsilon}_{\sigma_n \wedge t},Y^{\varepsilon}_{\sigma_n \wedge t})) &\leq& f(x_0,y_0) +  \E\int_0^{\sigma_n \wedge t} (2c_1 d \| X^{\varepsilon}_r\|  \|b \|_\infty + 2 c_2(M,R))  dr \nonumber\\ &&+  (m s(\varepsilon)^2/\varepsilon + d \varepsilon^{2\alpha}) \E(\sigma_n \wedge t) \nonumber\\
       &\leq & f(x_0,y_0) +  c_3\int_0^{t}  E\| X^{\varepsilon}_r\| dr  + c_4 t \nonumber\\
       &\leq & f(x_0,y_0) +  c_5\int_0^{t} (1+r) dr  + c_4 t \nonumber\\
       &\leq & c_6 + c_7t +c_8 t^2,\nonumber
       \end{eqnarray}
 where the penultimate inequality uses $c_3 \E[\|X^{\varepsilon}_r\|] \leq c_5(1+r)$ for a suitable $c_5$, a fact that follows from \eqref{ex} and the boundedness assumption on $b$ in (U1). As $\sigma_n \rightarrow \zeta$ almost surely, the above would imply
\begin{equation}
       \E[f(X^{\varepsilon}_{\zeta \wedge t},Y^{\varepsilon}_{\zeta \wedge t})] \leq  c_6 + c_7t +c_8 t^2.\end{equation}
Thus if $\zeta < t$ then we have a contradiction, as the left-hand side is infinity and the right-hand side is finite. As $t >0$ was arbitrary, we have $\zeta = \infty$ almost surely. This establishes nonexplosiveness of the process and completes the proof of strong uniqueness.
   \qed

\section{Nonlinear Filtering Equation}
\label{app:nonlinearfiltering}
Let {$x_0 \in\R^d$, $y_0\in\R^m$, $\sigma_1 >0$ and $\sigma_2
>0$. {On the probability space $(\Omega, {\cal F}, \P)$ let
  $\{B_t\}_{t \geq 0}$ and $\{W_t\}_{t \geq 0}$ be Brownian motions on
  $\R^d$ and $\R^m$ respectively.} In this section, we consider the
coupled diffusion $(X_t, Y_t)_{t \in [0,T]}$ on
$\mathbb{R}^d\times\mathbb{R}^m$ described by
\begin{eqnarray}
  X_t &=& x_0 + \int_0^t b_1(X_s, Y_s)ds + \sigma_1 B_t,
  \label{newxtimechange} \\
  Y_t &=& y_0 + \int_0^t b_2(X_s, Y_s)ds + \sigma_2 W_t,
  \label{newytimechange}
\end{eqnarray}
where $0 \leq t \leq T$, $b_1 : \mathbb{R}^d \times \mathbb{R}^m \rightarrow \mathbb{R}^d$ and $b_2 : \mathbb{R}^d \times \mathbb{R}^m \rightarrow \mathbb{R}^m.$ {We will make the following assumptions:}
\begin{itemize}
  \item $b_1 \in C_b(\R^d \times \R^m)$  is locally Lipschitz continuous in $y$-variable and is uniformly (w.r.t.\ $y$) Lipschitz continuous in $x$-variable, i.e.   $\exists K_1 > 0  $ such that $\forall \ x, x' \in \R^d, y \in \R^m$
 \begin{equation*}
    \norm{b_1(x, y) - b_1(x', y)}~ \leq K_1 \norm{x-x'}.
  \end{equation*}
\item  $b_2 \in C^1(\mathbb{R}^d \times \mathbb{R}^m).$ Further, $b_2(x, y)$ is  uniformly (w.r.t.\ $y$) Lipschitz continuous in $x$-variable, i.e.  $\exists K_2 > 0 $ such that $\forall \ x, x' \in \R^d$,  $ y \in \R^m$,
  \begin{equation*}
    \norm{b_2(x, y) - b_2(x', y)}~  \leq K_2 \norm{x-x'}.
  \end{equation*}
\item There exists $R>0, M>0$ such that for all $x \in \R^d$
    \begin{equation*}
 \sup_{\|y \| \leq R} \| b_2(x,y)\| \leq M.
    \end{equation*}
Further, there exist $K_4 >0$ and $R^\prime \geq R$ such that for all $x \in \R^d$
\begin{equation*}   - \langle b_2(x,y) , y \rangle > K_4 \|y\|^2, \quad \| y \| > R^\prime.\end{equation*}
\end{itemize}

Using the {above assumptions in the} same proof as in Lemma
\ref{lemma:eunexp}, it is standard to see that the above coupled
system has a unique strong solution. With $f \in C_b^2(\R^m)$, {for $y \in \R^m$ let}
\[
 { \LA^x_2 f  }(y) = \frac{\sigma_2^2}{2} \Delta f (y) + \langle b_2(x,\cdot), \nabla f (y) \rangle,
\]
where $x$ is treated as a parameter. Let $\FA_t^X = \sigma(X_s, s \leq t)$, $\FA^X = \bigvee_{t \geq 0} \FA_t^X$. Define $\FA_t^{X,Y}$ and $\FA^{X,Y}$ analogously. Define $\pi_t(dy)$ as the conditional law of $Y_t$ given $\FA_t^X$ so that
\[
\pi_t(f) := \ E[f(Y_t) | \FA_t^X] \quad \mbox{ for } f \in C_b^2(\R^m).
\]

\begin{proposition}\label{fkksde} {\bf (Nonlinear Filtering Equation)}
The measure valued process $\pi$ is the unique solution to the (Fujisaki-Kallianpur-Kunita) nonlinear filtering equation
\begin{equation}\label{KushnerStratonovich}
 \pi_t(f) \ = f(y_0) + \ \int_0^t \pi_s({\LA^{X_s}_2  f)} ds
+ {\frac{1}{\sigma_1}} \int_0^t\langle \pi_s(f b_1(X_s, \cdot)) -
\pi_s (f) \pi_s(b_1(X_s, \cdot)), d \tilde{B}_s \rangle, \quad f \in C_b^2(\R^m),
\end{equation}
where $\tilde{B}$ is a standard Brownian motion.
\end{proposition}

\begin{remark}
${\tilde{B}}$ is  explicitly defined later in the proof. It is called the `innovations process' and, under mild technical conditions, is known to generate the same increasing $\sigma$-fields as $B$ \cite{AM81}.
\end{remark}

 We will closely mimic the arguments in [\cite{BC09}, Chapter 3]
 proved for the case when $b_1(x,y)$ is a function of the first
 argument alone. In our setting, $b_1(x,y)$ is a function of both
 arguments.

Set
\begin{equation}\label{Lambda}
\Lambda_s \ = \ \exp \left\{-\frac{1}{\sigma_1} \int^s_0 \langle b_1(X_u , Y_u), d B_u \rangle
  - \frac{1}{2\sigma_1^2} \int^s_0 \norm{b_1(X_u, Y_u)}^2 du
\right\}, \ s \geq 0,
\end{equation}
and
\[
\Bar{B}_s \ = \ B_s + \frac{1}{\sigma_1} \int^s_0 b_1(X_u , Y_u) du
, \ s \geq 0.
\]
Define the probability measure $Q$ by
\[
\frac{d Q\Big|_{\FA_s^{X,Y}}}{dP\Big|_{\FA_s^{X,Y}}} \ = \ \Lambda_s, s > 0.
\]
This consistently defines $Q$ on $\FA^{X,Y}$. {As $b_1$ is bounded,} by the
Cameron-Martin-Girsanov theorem, it follows that $\Bar{B}_\cdot$ is
an $\mathbb{R}^d$-valued standard {Brownian motion} under
$Q$.  Under $Q$, the joint process $(X, Y)$ given by
(\ref{newxtimechange}) - (\ref{newytimechange}) takes the form
{\begin{eqnarray}\label{mainsdesundernewavthar}
 X_t & = & x_0 + \sigma_1  \Bar{B}_t, \nonumber \\
 Y_t & = & y_0 + \int_0^tb_2(X_s, Y_s) ds + \sigma_2  W_t.
\end{eqnarray}
}

{Before we begin the proof we need some preliminary lemmas.}

\begin{lemma}\label{ap1} For $t > 0$, let $Z$ be a $Q$-integrable $\FA_t^{X,Y}$-measurable
$\mathbb{R}^d$-valued random variable. Then
\[
E^Q [Z | \FA_t^X] \ = \ E^Q [Z | \FA^X].
\]
\end{lemma}

{\em Proof:} Set
\[
\tilde{\mathcal F}^X_t \ = \ \sigma(X_{t+s}- X_t, s \geq 0).
\]
Then ${\mathcal F}^X = \tilde{\mathcal F}^X_t\vee{\mathcal F}^X_t$, and
since $X_s = \sigma_1 \Bar{B}_s$, an $\{{\mathcal
  F}_s^X\}$-Wiener process under $Q$, $\tilde{\mathcal F}^X_t$ is
independent of $\FA_t^X$ under $Q$. Hence
\[
\begin{array}{lll}
E^Q [ Z | \FA_t^X ] & = & E^Q [ Z | \tilde{\mathcal F}^X_t\vee
  {\mathcal F}_t^X ] \\ & = & E^Q[ Z | {\mathcal F}^X] . \\
\end{array}
\]
This completes the proof of the lemma.  \qed

\begin{lemma}\label{ap2} Let $\{{\alpha}_t, ~t \geq 0\}$ be an $\{\FA_t^{X,Y}\}$-progressively measurable $\mathbb{R}$-valued process such that
\[
E^Q \Big[ \int^t_0 {\alpha}^2_s ds \Big] < \infty \ \forall \ t > 0.
\]
Then
\[
E^Q \Big[ \int^t_0 {\alpha}_s dX_s \Big| {\mathcal F}^X \Big] \ =
\ \int^t_0 E^Q[ {\alpha}_s | {\mathcal F}^X] dX_s .
\]
\end{lemma}

{\em Proof:} Using Lemma \ref{ap1}, it follows that
\[
E^Q \Big[ \int^t_0 {\alpha}_s dX_s \Big| {\mathcal F}^X \Big], \ E^Q
[{\alpha}_t | {\mathcal F}^X]
\]
are ${\mathcal F}^X_t$-measurable. Hence using the `density result' of
Krylov and Rozovskii, see [\cite{BC09}, Lemma B.39, p.355], it is
enough to show
\begin{eqnarray}
E^Q \Big[ \beta_t E^Q \Big[ \int^t_0 {\alpha}_s dX_s \Big| {\mathcal
      F}^X \Big] \Big] & = & E^Q \Big[ \beta_t \int^t_0 E^Q [{\alpha}_s
    | {\mathcal F}^X] dX_s \Big] \\ \nonumber
\end{eqnarray}
for all process $\beta(\cdot)$ of the form
\[
\beta_t = 1 + \int^t_0 i \langle \beta_s r_s, dX_s \rangle
\]
for a deterministic $r \in L^{\infty}([0, t] ;
\mathbb{R}^d)$. Consider
\begin{align*}
 E^Q \Big[ \beta_t E^Q \Big[ \int^t_0 {\alpha}_s dX_s
     \Big| {\mathcal F}^X \Big] \Big]  & =  E^Q
  \Big[ \beta_t \int^t_0 {\alpha}_s dX_s \Big] \\ & =
  E^Q \Big[ \int^t_0 {\alpha}_s dX_s \Big] + E^Q \Big[ \Big( \int^t_0 i
    \langle \beta_s r_s, dX_s \rangle \Big) \Big( \int^t_0 {\alpha}_s
    dX_s \Big) \Big]\\
  &=  \sigma_1^2 E^Q
  \Big[\int^t_0i \beta_s r_s {\alpha}_s ds \Big]\\&=
  \sigma_1^2 E^Q \Big[E^Q \Big[\int^t_0i \beta_s r_s
      {\alpha}_s ds \Big| {\mathcal F}^X \Big] \Big]
  \\ &=  E^Q \Big[ \Big( \int^t_0 i \langle \beta_s r_s, dX_s \rangle \Big)
    \Big( \int^t_0 E^Q [{\alpha}_s | {\mathcal F}^X] dX_s \Big)
    \Big]\\ &=  E^Q \Big[ \beta_t \int^t_0 E^Q
    [{\alpha}_s | {\mathcal F}^X] dX_s \Big] \\
\end{align*}
This completes the proof of the lemma.  \qed

\begin{lemma}\label{ap3} Let $x \in \R^d$. Let $\{{\alpha}_t, ~t \geq 0\}$ be $\{\FA_t^{X,Y}\}$-progressively
measurable process such that
\[
E^Q \Big[ \int^t_0 {\alpha}^2_s d \langle M^f \rangle_s \Big] < \infty,
\ f \in C^2_b(\mathbb{R}^m), \ t \geq 0,
\]
where
\[
M^f_t \, = \, f(Y_t) - f({y_0}) - \int^t_0
\LA_2^x f(Y_s) ds,
\]
and $\langle M^f \rangle_t$ is its quadratic
variation.  Then
\[
E^Q \Big[ \int^t_0 {\alpha}_s d M^f_s \Big| {\mathcal F}^X \Big] = 0.
\]
\end{lemma}

{\em Proof:} Via It\^{o}'s formula, we first obtain
\begin{equation}
\label{eqn:dM}
dM_t^f = \sigma_2 \langle \nabla f(Y_t), dW_t \rangle.
\end{equation}
Under $P$, this is driven by a Brownian motion independent of $B_t$, which leads to $$\langle M^f, X \rangle_t = 0, \ P-\mbox{ almost surely}$$ and hence $Q-$almost surely. Using this, the proof follows along the lines of the proof of Lemma \ref{ap2}.  \qed

Set
\[
\tilde{\Lambda}_t = \Lambda^{-1}_t, \ t \geq 0,
\]
{ and for $g \in C^2(\mathbb{R}^d \times \mathbb{R}^m)$ with a little abuse of notation denote $$ \pi_t(g) := \pi_t(g(X_t, \cdot)) = E [ g(X_t, Y_t) | {\mathcal F}_t^X]$$}
We then have the following.

\begin{lemma}\label{ap4} (Kallianpur-Striebel formula)
For $g \in C^2(\mathbb{R}^d \times \mathbb{R}^m)$,
\[
\pi_t(g)  = \ \frac{E^Q [ \tilde{\Lambda}_t g | {\mathcal
      F}^X]}{E^Q[\tilde{\Lambda}_t | {\mathcal F}^X]}.
\]
\end{lemma}

{\em Proof:}  In view
of Lemma \ref{ap1}, it is enough to show that
\[
\pi_t(g) E^Q[\tilde{\Lambda}_t | \FA_t^X] \ =
\ E^Q[\tilde{\Lambda}_t g | \FA_t^X].
\]
Since both left and right sides are $\FA_t^X$-measurable, it is enough to show that

\begin{equation}\label{lemma4.4eq1}
E^Q \left[ \beta \pi_t(g) E^Q \left[\tilde{\Lambda}_t | \FA_t^X \right] \right] \ =
\ E^Q \left[ \beta \tilde{\Lambda}_t g \right]
\end{equation}
for all $\FA_t^X$-measurable $\beta$. This is now easily verified since, for such $\beta$, we have
\begin{eqnarray*}
  E^Q \left[ \beta \pi_t(g) E^Q \left[ \tilde{\Lambda}_t | \FA_t^X \right] \right] &=&  E^Q \left[ \beta \pi_t(g) \tilde{\Lambda}_t \right]
   =  E \left[ \beta \pi_t(g)  \right]  = E \left[ \beta E \left[ g | \FA_t^X \right]  \right]   =  E \left[ E \left[ \beta g | \FA_t^X \right]  \right] \\
  & = & E \left[ \beta g \right]    =  E^Q \left[ \tilde{\Lambda}_t \beta g  \right].
\end{eqnarray*}
This completes the proof of the lemma.  \qed

{We are now ready to prove Proposition \ref{fkksde}.}
We shall derive first the Zakai equation solved by certain unnormalized conditional laws. Then we shall show existence to the Fujisaki-Kallianpur-Kunita) nonlinear filtering equation \eqref{KushnerStratonovich}, followed by uniqueness.

{\em Proof of Proposition \ref{fkksde}:} Observe that $\{\Lambda_t, t \geq 0\}$ is given by the solution of the
SDE
\[
{\Lambda_t \ =\ 1 - \int_0^t\Lambda_s \sigma_1^{-1} \langle b_1(X_s,Y_s), dB_s \rangle,}
\]
for $t \geq 0$.
Hence by a routine application of It$\hat{\rm o}$'s formula it follows
that
\begin{equation}\label{repLambda1}
{\tilde{\Lambda}_t \ = \ 1 + \int_0^t\tilde{\Lambda}_s \sigma_1^{-2}
\langle b_1(X_s, Y_s), dX_s \rangle.}
\end{equation}
From this, since $X_t$ is driven by $B_t$ and $Y_t$ is driven by $W_t$, for $f \in C_b^2(\R^m)$, the cross-variation $\langle \tilde{\Lambda}, f(Y_{\cdot}) \rangle_t = 0$ $P$-a.s. and hence $Q$-a.s.  Using It$\hat{\rm o}$'s formula again, we get
 \[
{\tilde{\Lambda}_t f(Y_t) = f(y_0) + \int_0^t\tilde{\Lambda}_s [
{ \LA^{X_s}_2  f (Y_s) }ds +
  \sigma_2 \langle \nabla f(Y_s), dW_s \rangle ] + \int_0^t f(Y_s) d \tilde{\Lambda}_s ,}
\]
and hence, using (\ref{eqn:dM}) and (\ref{repLambda1}), we get
\begin{equation}\label{tleq}
\tilde{\Lambda}_t f(Y_t) = f(y_0) +
  \int^t_0 \tilde{\Lambda}_s {\LA^{X_s}_2  f (Y_s) } ds + \int^t_0 \tilde{\Lambda}_s d M^f_s
 + \sigma_1^{-2} \int^t_0 \tilde{\Lambda}_s f(Y_s)
  \langle b_1(X_s, Y_s), dX_s \rangle.
  \end{equation}
Taking conditional expectation $E^Q[ \ \cdot \ | {\mathcal F}^X]$ in (\ref{tleq}) we have using Lemma \ref{ap3}  we  have
\begin{align}
E^Q\left [ \tilde{\Lambda}_t f(Y_t)| {\mathcal F}^X \right ]   &= f(y_0) + E^Q\left [ \int^t_0
\tilde{\Lambda}_s( {\LA^{X_s}_2  f(Y_s) )  } ds| {\mathcal F}^X \right ] + \sigma_1^{-2}E^Q\left [
  \int^t_0 \tilde{\Lambda}_sf((Y_s) \langle   b_1(X_s, Y_s)), d X_s \rangle  | {\mathcal F}^X \right ], \nonumber\\
&\mbox{ and using Lemma \ref{ap2} we have the above is}\nonumber\\
&=f(y_0) +  \int^t_0
E^Q\left [ \tilde{\Lambda}_s( {\LA^{X_s}_2  f(Y_s) )  } | {\mathcal F}^X \right ]ds + \sigma_1^{-2}
\int^t_0 \langle E^Q\left [ \tilde{\Lambda}_sf((Y_s)    b_1(X_s, Y_s))| {\mathcal F}^X \right ], d X_s \rangle  ,
\end{align}
{  For $g \in C(\R^d\times \R^m)$ denoting
\[
\rho_t(g) \ = \ \pi_t(g) E^Q[ \tilde{\Lambda}_t | {{\mathcal F}^X} ].
\]}
in (\ref{tleq}) and using Lemma \ref{ap4} we arrive at the Zakai equation
\begin{equation} \label{zeq} \rho_t(f) = f(y_0) + \int^t_0
\rho_s( {\LA^{X_s}_2  f )  } ds + \sigma_1^{-2}
\int^t_0 \langle \rho_s( f b_1(X_s, \cdot)), d X_s \rangle.
\end{equation}
For $\1 :=$ the constant function identically equal to
$1$, we see that $\rho_t(\1) = E^Q[ \tilde{\Lambda}_t | {{\mathcal F}^X} ]$, and hence
 \begin{equation}
\pi_t(f) = \frac{\rho_t(f)}{\rho_t(\1)}. \label{unn}
\end{equation}
The nonnegative measure valued process $\{\rho_t\}_{t \geq 0}$ is called the process of
unnormalized conditional laws in view of (\ref{unn}).

Now we are ready to prove the
existence theorem for the Fujisaki-Kallianpur-Kunita (FKK) equation,
(\ref{KushnerStratonovich}). From the Zakai equation (\ref{zeq}) we
get
\begin{equation}\label{Theorem4.2eq1}
\rho_t(f) \, = \, f(y_0) + \int^t_0
\rho_s(\1)\pi_s( {\LA^{X_s}_2  f )  } ds + \sigma_1^{-2}
\int^t_0 {\rho_s(\1)}\langle \pi_s( f b_1(X_s, \cdot)), d X_s \rangle ,
\end{equation}
In particular, one can deduce that
\begin{equation}\label{Theorem4.2eq2n}
\rho_t(\1)  = 1 + \sigma_1^{-2} \int_0^t   \rho_s(\1)
\langle \pi_s(  b_1(X_s, \cdot)), dX_s \rangle.
\end{equation}
Using It$\hat{\rm o}$'s formula, we get
\begin{equation} \label{Theorem4.2eq2} \frac{1}{\rho_t(\1)} = 1 - \sigma_1^{-2}\int_0^t\frac{1}{\rho_s(\1)}
\langle \pi_s(  b_1(X_s, \cdot)), dX_s \rangle +  \sigma_1^{-2}\int_0^t\frac{1}{\rho_s(\1)} \| \pi_s(b_1(X_s,\cdot))\|^2 ds,\end{equation}
Note that the cross-variation
\begin{equation}\label{Theorem4.2eq4}
 \langle \rho(f) , \frac{1}{\rho(\1)} \rangle_t \, = \, -\int_0^t
\sigma_1^{-2} \langle \pi_s( b_1), \pi_s( b_1 f) \rangle ds,
\end{equation}
 It$\hat{\rm o}$'s formula, for the product of $\rho_t(f)$ and $\frac{1}{\rho_t(\1)}$ we get
  \begin{align}\nonumber
 \frac{\rho_t(f)}{\rho_t(\1)}  &=  f(y_0) + \int_0^t\rho_s(f) d \frac{1}{\rho_s(\1)} + \int_0^t \frac{1}{\rho_s(\1)} d\rho_s(f)  + \langle \rho(f), \frac{1}{\rho(\1)}\rangle_t
\end{align}
  Substituting (\ref{Theorem4.2eq1}),(\ref{Theorem4.2eq2}), and (\ref{Theorem4.2eq4}) in the above we have
   \begin{align}
 \frac{\rho_t(f)}{\rho_t(\1)}  &=  f(y_0) + \int_0^t\rho_s(f) \left[-\frac{1}{\rho_s(\1)}
\sigma_1^{-2}\langle \pi_s(  b_1(X_s, \cdot)), dX_s \rangle + \sigma_1^{-2}\frac{1}{\rho_s(\1)} \| \pi_s(b_1(X_s,\cdot))\|^2 ds, \right]  \nonumber \\ &\hspace{1in}+ \int_0^t \frac{1}{\rho_s(\1)}\left [ \rho_s(\1)\pi_s( {\LA^{X_s}_2  f )  } ds + \sigma_1^{-2}{\rho_s(\1)}\langle \pi_s( f b_1(X_s, \cdot)), d X_s \rangle \right] \nonumber \\&   \hspace{1in}-    \int_0^t
 \sigma_1^{-2} \langle \pi_s( b_1), \pi_s( b_1 f) \rangle ds. \nonumber
   \end{align}
From (\ref{unn}) and simple algebra in the above we have
   \begin{align}  \label{Theorem4.2eq3}
\pi_t(f) &=f(y_0) + \int_0^t \pi_s( {\LA^{X_s}_2  f )  } ds + \sigma_1^{-2}\int_0^t \langle \pi_s( f b_1)- \pi_s(f) \pi_s(  b_1), d X_s -\pi_s(b_1)ds\rangle
\end{align}
Let
\[
I_t = X_t - \int_0^t\pi_s( b_1) ds,
\]
the so called `innovation process'. For $0 \leq s < t $, we have
\[
\begin{array}{lll}
E [ I_t - I_s | \mathcal{F}_s^X] & = & \displaystyle{ E \Big[ \int^t_s
    E[ b_1(X_u, Y_u) - \pi_u( b_1(X_u, \cdot)) \ | \ {\mathcal F}_u^X ] du
    \Big| \mathcal{F}_s^X\Big]}\\ & = & \displaystyle{ \int^t_s E \Big[
    b_1(X_u, Y_u) - E [ b_1(X_u, Y_u) \ | \ {\mathcal F}_u^X]| \mathcal{F}_s^X
    \Big] du} \\ & = & 0.\\
\end{array}
\]
Thus $\{{I}_t | t \geq 0\}$ is an $\{{\mathcal F}_t^X\}$-martingale
with mean $0$ and quadratic variation $\sigma_1^2t$. Thus by Levy's
characterization, $I$ is a scaled Brownian motion. Define
\begin{equation}
\tilde{B}_t := \sigma_1^{-1}I_t, \ t \geq
0. \label{newinnovation}
\end{equation}
So $\tilde{B}_t$ is a $\{{\mathcal  F}_t^X\}$-adapted standard Brownian motion under $P$.
Therefore we have shown that,
   \begin{align}  \label{Theorem4.2eq10}
\pi_t(f) &=f(y_0) + \int_0^t \pi_s( {\LA^{X_s}_2  f )  } ds + {\sigma_1^{-1}}\int_0^t \langle \pi_s( f b_1)- \pi_s(f) \pi_s(  b_1), d\tilde{B}_s\rangle,
\end{align}
with $\tilde{B}_s$ being a standard Brownian motion. Thus we have shown existence of a solution to  the FKK equation. Uniqueness of the FKK equation in the sense of martingale problem follows from Theorem 3.3 of Kurtz and Ocone \cite{KO88}. Note that while Kurtz and Ocone \cite{KO88} cite nonlinear filtering as an example of this theorem, they consider the classical formulation (see \cite[Theorem 4.1]{KO88}) which is more restrictive than ours. However  the aforementioned theorem (\cite[Theorem 3.3]{KO88}) is general enough to cover our problem.

 \qed

{From (\ref{Theorem4.2eq2n}) we have
$$\rho_t(\1) = \exp \left\{\sigma_1^{-{2}}\int_0^t\langle\pi_s(b_1(X_s,
  \cdot )), dX_s\rangle -
  \frac{\sigma_1^{-{4}}}{2}\int_0^t\|\pi_s(b_1(X_s, \cdot
  ))\|^2ds \right\}.
$$
Since $\pi_t(f) = \frac{\rho_t(f)}{\rho_t(\1)},$
$$\rho_t(f) = \pi_t(f)\rho_t(\1) =
\pi_t(f) \exp \left\{\sigma_1^{-{2}}\int_0^t\langle\pi_s(b_1(X_s, \cdot )),
  dX_s\rangle - \frac{\sigma_1^{-{4}}}{2}\int_0^t\|\pi_s(b_1(X_s,
  \cdot ))\|^2ds\right\}.$$
Thus solutions $\pi, \rho$ of FKK, resp.\ Zakai
equations are in one-one correspondence and uniqueness of one implies
that of the other.}
\begin{remark} \label{rem:filtering}
It is interesting to note that some of the earlier uniqueness arguments for the classical framework such as one using multiple Wiener integral expansion due to \cite{KU81} or via the Clark-Davis `pathwise' filter as in \cite{H85}, do not work for our case. (The latter would work only if $b_1(x, \cdot ) = \nabla F(x, \cdot)$ for a suitable $F$.)
\end{remark}

\section{Laplace's principle}
\label{appendix:Laplace}

We now characterize weak limit points of the sequence of invariant measures for the fast process $\nu^{\varepsilon_n, x^{\varepsilon_n}}$ when $\varepsilon_n \to 0$ and for deterministic $x^{\varepsilon_n} \to x$. This is used in the proof of Proposition \ref{characterisation}(c,d).

\begin{lemma} \label{l:ieu}
  Let $n \geq 1, 0 < \varepsilon_n < 1, x^{\varepsilon_n} \in {\mathbb R}^d$, and $x \in \mathbb{R}^d$. Suppose $\varepsilon_n \rightarrow 0$ and $x^{\varepsilon_n} \rightarrow x $ as $n\rightarrow \infty$.
  \begin{enumerate}
    \item[(a)] Then  the sequence of measures $\nu^{\varepsilon_n, x^{\varepsilon_n}}$ is tight and any limit point is  supported on $\arg\min\{U(x, \cdot)\}$.
 \item[(b)] Assume (U4) and let $x \in D_L^\circ$ for some $L \geq 1$. Then $\nu^{\varepsilon_n, x^{\varepsilon_n}}$ converges weakly to $\nu^{0,x}$,  where $\nu^{0,x}$ is  given by
       $$\sum_{i=1}^{L} \delta_{y_i(x)}\frac{\left(\mbox{Det}\left [D_y^2 U(x,y_i(x))\right]\right)^{-\frac{1}{2}}}{\sum_{j=1}^{L} \left(\mbox{Det}\left [D_y^2 U(x,y_j(x))\right] \right)^{-\frac{1}{2}}}.$$
\end{enumerate}
  \end{lemma}
\vfill
    {\em Proof :} (a) Using (\ref{U21}) and (\ref{U22}) it is easy to see (see Lemma \ref{l:U24} in Appendix A) that there is $K_4 >0$ {and an $R^\prime \geq R$ such that
\begin{equation}   \langle \nabla_y U(x, y) , y \rangle > K_4 \|y\|^2, \quad \|y\| > R^\prime.\label{U24a}\end{equation}
}
Let  $h: \R^m \rightarrow \R$ be given by $$ h(y) ~ = ~ \parallel y \parallel^2.$$
 Using (\ref{U24a}), there is $R^\prime >0$ such that
\bena
 &&{\LA^{\varepsilon_n, x^{\varepsilon_n}}( h )} (y)= \frac{s(\varepsilon_n)^{2}}{2}\Delta h(y) - \langle \nabla_yU(x^{\varepsilon_n}, y), \nabla h(y) \rangle = m s(\varepsilon_n)^{2}  - 2 \langle \nabla_yU(x^{\varepsilon_n}, y), y \rangle \nonumber\\
 && < m s(\varepsilon_n)^{2}  - 2K_4 \parallel y \parallel^2 < 0 , \label{h:lyap}
\eena
for all  $\parallel y \parallel > R^\prime$ and $n \geq 1.$ We can assume without loss of generality that $\max\{\| x^{\varepsilon_n} \|, \|x\|\} \leq R^\prime$. So by regularity assumption on $U$ from (U1) we have

\begin{equation} \int_{\parallel y \parallel  \leq R^\prime } \mid  m s(\varepsilon_n)^{2}  - 2 \langle \nabla_yU(x^{\varepsilon_n}, y), y \rangle \mid \nu^{\varepsilon_n, x^{\varepsilon_n}}(dy) \leq  K
 \label{h:mprhyp} \end{equation}
for some $K \equiv K(m,s,R^\prime,U)>0$. Using (\ref{h:lyap}), (\ref{h:mprhyp}) along with  Proposition \ref{characterisation} in \cite{MPR05} and its proof, we have for all $n \geq 1$
\begin{equation} \label{h:L1} \nu^{\varepsilon_n, x^{\varepsilon_n}}( \mid \LA^{\varepsilon_n}(x^{\varepsilon_n} , h ) \mid)  < 2K.  \end{equation}

Define $g: \R^m \rightarrow \R$ by $g(y) = 2K_4 \parallel y \parallel^2 - m s(\varepsilon)^{2}$.  Using (\ref{h:lyap}) and (\ref{h:L1}) we have
\begin{equation}\label{nutight}
  \int_{\parallel y \parallel > R^\prime} g(y) \nu^{\varepsilon_n, x^{\varepsilon_n}}(dy) \leq   \nu^{\varepsilon_n, x^{\varepsilon_n}}( \mid \LA^{\varepsilon_n}(x^{\varepsilon_n} , h ) \mid)  < 2K.
\end{equation}
As $g(y) \rightarrow \infty$ when $\parallel y \parallel \rightarrow
\infty$ we can conclude that the sequence of measures
$\{\nu^{\varepsilon_n, x_n}\}_{n \geq 1}$ is tight. We will now show that any limit point $\nu$ is supported on $\arg\min U(x, \cdot)$.

Let $z \in \R^m, z \not \in \arg \min \{U(x, \cdot) \}.$ As $U(x^{\varepsilon_n},\cdot)$ converges to $U(x,\cdot)$ uniformly on compact sets, there exists $\delta >0 $ and $r >0$ such that
\[ U(x^{\varepsilon_n},y) > U(x,y_i(x)) + \frac{\delta}{2}, \,\, \forall y \in B(z,r)\]
and
\[ U(x^{\varepsilon_n},y) < U(x,y_i(x)) + \frac{\delta}{4}, \,\, \forall y \in B(y_i(x),r). \]
Therefore, for $n \geq 1$,
\begin{eqnarray*}
  \frac{\nu^{\varepsilon_{n}, x^{\varepsilon_n}}(B(z,r))}{\nu^{\varepsilon_{n}, x^{\varepsilon_n}}(B(y_i(x), r))} &=& \frac{\int_{B(z,r)} e^{-2\frac{U(x^{\varepsilon_n},y)}{s(\varepsilon_n)^2}} dy}{\int_{B(y_i(x),r)} e^{-2\frac{U(x^{\varepsilon_n},y)}{s(\varepsilon_n)^2}}dy} =  \frac{\int_{B(z,r)} e^{-2\frac{U(x^{\varepsilon_n},y) -U(x,y_i(x))}{s(\varepsilon_n)^2}} dy}{\int_{B(y_i(x),r)} e^{-2\frac{U(x^{\varepsilon_n},y) -U(x,y_i(x))}{s(\varepsilon_n)^2}}dy}\\
  &\leq & \frac{\mid B(z,r)\mid  e^{-\frac{\delta}{ s(\varepsilon_n)^2}}}{\mid B(y_i(x),r)\mid e^{-\frac{\delta}{2s(\varepsilon_n)^2}}} = e^{-\frac{\delta}{2 s(\varepsilon_n)^2}}.
\end{eqnarray*}
Therefore, $$\displaystyle \lim_{n \rightarrow \infty} \nu^{\varepsilon_{n}, x^{\varepsilon_n}}(B(z,r)) = 0.$$

Hence any limit point $\nu$  is supported on the $\arg \min \{ U(x,\cdot)\}.$

(b) Let $x \in D_L^\circ$ for some $L \geq 1$. Under (U4) the global
minima $y_i(x), 1 \leq i \leq L$ are nondegenerate, i.e., the matrix
$D^2_yU(x, y_i(x))$ is positive definite for $1 \leq i \leq
L$. Since $D_L^\circ$ is open, using $U \in C^2(\R^m \times \R^d)$ in (U1), with a suitable relabelling of the minima if necessary, we have $y_i(x^{\varepsilon_n}) \to y_i(x) \ \forall 1 \leq i \leq L$ as $n \to \infty$. Let
$B_i$ be the ball centered at $y_i(x)$ with radius $1$ for each $i$. Let $n$ be sufficiently large so that $y_i(x^{\varepsilon_n})$ are in a  ball centered at $y_i(x)$ with radius $\frac{1}{2}$.  Let $B^n_i$ be ball centered at $y_i(x^{\varepsilon_n})$ with radius $\frac{1}{4}$ for each $i$.
Note that $\nabla_yU(x^{\varepsilon_n}, y_i(x^{\varepsilon_n})) = 0$  and  $U(x^{\varepsilon_n}, y_i(x^{\varepsilon_n})) = \min U(x^{\varepsilon_n}, \cdot):= u_{\min}$ (say) as it does not depend on $i$. Using Taylor's
expansion up to second order, we have that for each $ y \in B_i$ there is a $\tilde{y}_i(x^{\varepsilon_n}) \in B_i$ such that
\begin{eqnarray*}
  U(x^{\varepsilon_n}, y) &=& U(x^{\varepsilon_n}, y_i(x^{\varepsilon_n})) + \frac{1}{2}(y - y_i(x^{\varepsilon_n}))^TD^2_yU(x^{\varepsilon_n}, \tilde{y}_i(x^{\varepsilon_n}))(y - y_i(x^{\varepsilon_n}))\\
  &=& u_{\min} + \frac{1}{2} (y - y_i(x^{\varepsilon_n}))^TD^2_yU(x^{\varepsilon_n}, \tilde{y}_i(x^{\varepsilon_n}))(y - y_i(x^{\varepsilon_n})).
\end{eqnarray*}

The above and  standard fact about Gaussian random variables implies:
\begin{align}
  \int_{B_i}e^{-2\frac{U(x^{\varepsilon_n}, y)}{s(\varepsilon_n)^2}}dy &\leq e^{-2\frac{u_{\min}}{s(\varepsilon_n)^2}} \int_{\R^m}e^{-2\frac{(y - y_i(x^{\varepsilon_n}))^TD^2_yU(x^{\varepsilon_n}, \tilde{y}_i(x^{\varepsilon_n}))(y - y_i(x^{\varepsilon_n}))}{2s(\varepsilon_n)^2}}dy \nonumber\\& = e^{-2\frac{u_{\min}}{s(\varepsilon_n)^2}} \left( \frac{(2\pi)^m s(\varepsilon_n)^2}{2}\mbox{Det}\left(D^2_yU(x^{\varepsilon_n}, y_i(x^{\varepsilon_n}))^{-1}\right)\right)^{\frac{1}{2}}; \label{ieubs1}\end{align}
 and
\begin{align}
  \int_{B_i}e^{-2\frac{U(x^{\varepsilon_n}, y)}{s(\varepsilon_n)^2}}dy &\geq e^{-2\frac{u_{\min}}{s(\varepsilon_n)^2}} \int_{B^n_i}e^{-2\frac{(y - y_i(x^{\varepsilon_n}))^TD^2_yU(x^{\varepsilon_n}, \tilde{y}_i(x^{\varepsilon_n}))(y - y_i(x^{\varepsilon_n}))}{2s(\varepsilon_n)^2}}dy\nonumber \\& = e^{-2\frac{u_{\min}}{s(\varepsilon_n)^2}} \left(\frac{(2\pi)^ms(\varepsilon_n)^2}{2}\mbox{Det}\left(D^2_yU(x^{\varepsilon_n}, y_i(x^{\varepsilon_n}))^{-1}\right)\right)^{\frac{1}{2}} P\left( Z^m \in \frac{1}{s(\varepsilon_n)} A_{i,n}\right),\label{ieubs2}\end{align}
where $A_{i,n} = \left\{ z \in \R^m : \left\| \left(D^2_yU(x^{\varepsilon_n}, y_i(x^{\varepsilon_n})) \right) ^{-1/2} z \right\| \leq \frac{1}{2\sqrt{2}} \right\}$ and $Z^m $ is a standard $m-$dimensional Gaussian random variable.
{ Note that $A_{i,n}$ is a bounded set in $\R^m$ and  by (U1), $U \in C^2(\R^d \times \R^m)$. So  as $n \rightarrow \infty$ we have  }
$$P \left( Z^m \in \frac{1}{s(\varepsilon_n)} A_{i,n} \right) \rightarrow 1 \mbox{ and } \left(\mbox{Det}\left(D^2_yU(x^{\varepsilon_n}, y_i(x^{\varepsilon_n}))^{-1}\right)\right)^{\frac{1}{2}} \rightarrow \left(\mbox{Det}\left(D^2_yU(x, y_i(x))^{-1}\right)\right)^{\frac{1}{2}}$$ for all $i$.
{Therefore using a standard sandwich argument we can conclude that, for balls $B_i$ and $B_j$, $1 \leq i, j \leq L$,
\begin{equation}
\label{eqn:Laplace3}
\frac{\nu^{\varepsilon_n,x^{\varepsilon_n}}(B_i)}{\nu^{\varepsilon_n,x^{\varepsilon_n}}(B_j)}= \frac{\int_{B_i}e^{-2\frac{U(x^{\varepsilon_n}, y)}{s(\varepsilon_n)^2}}dy}{\int_{B_j}e^{-2\frac{U(x^{\varepsilon_n}, y)}{s(\varepsilon_n)^2}}dy} \to \frac{\left(\mbox{Det}\left(D^2_yU(x, y_i(x))^{-1}\right)\right)^{\frac{1}{2}}}{\left(\mbox{Det}\left(D^2_yU(x, y_j(x))^{-1}\right)\right)^{\frac{1}{2}}} \,\, \mbox{ as } n \rightarrow \infty.
\end{equation}
From (a) we know that the sequence of measures $\{\nu^{\varepsilon_n,x^{\varepsilon_n}}\}_{n \geq 1}$ are tight and all limit points are measures supported on the $\arg\min U(x, \cdot)$. Consequently by (\ref{eqn:Laplace3}) we have that any limit point $\nu^{0,x}$ is given by
    $$\nu^{0,x}(\cdot) = \sum_{i=1}^{L}\frac{\left(\mbox{Det}\left
  [D_y^2 U(x,y_i(x))\right]\right)^{-\frac{1}{2}}}{\sum_{j=1}^{L}
  \left(\mbox{Det}\left [D_y^2 U(x,y_j(x))\right]
  \right)^{-\frac{1}{2}}}\delta_{y_i(x)}( \cdot).$$
Since all subsequential limit points are the same we have the result.}
\qed

}
\end{appendices}

{\bf Acknowledgements:} Research of S.R.A. was supported in part by ISF-UGC grant, research of V.S.B. was supported in part by a J.\ C.\ Bose Fellowship, research of K.S.K. was supported in part by the grant MTR/2017/000416 from SERB and research  of R.S. was supported in part by RBCCPS-IISc. S.R.A., V.S.B. and R.S. would like to thank the International Centre for Theoretical Sciences (ICTS) for hospitality during the  {\em Large deviation theory in statistical physics:  Recent advances and future challenges (Code:ICTS/Prog-ldt/2017/8).} The authors thank  Sanjoy Mitter for pointing out the reference \cite{KO88}, Laurent Miclo and  Patrick Cattiaux for suggestions on the spectral gap estimate in Proposition \ref{p:spgp}(b), and Konstantinos Spiliopoulos for pointing out several references in the literature.


 \noindent {\bf Siva Athreya}\\
Stat-Math Unit, Indian Statistical Institute, 8th Mile, Mysore Road,
     Bangalore 560059, India.
     Email: \texttt{athreya@isibang.ac.in}

     \medskip

\noindent{\bf Vivek S. Borkar}\\
   Department of Electrical Engineering, Indian Institute of Technology, Powai, Mumbai 400076, India.
     Email: \texttt{borkar.vs@gmail.com}

     \medskip

\noindent{\bf K. Suresh Kumar}\\
Department of Mathematics, Indian Institute of Technology, Powai, Mumbai 400076, India.
Email: \texttt{suresh@math.iitb.ac.in}

\medskip

\noindent{\bf Rajesh Sundaresan}\\
Department of Electrical Communication Engineering and Robert Bosch Centre for Cyber-Physical Systems, Indian Institute of Science, Bangalore 560012, India.
Email: \texttt{rajeshs@iisc.ac.in}

\end{document}